\documentclass{amsart}
\usepackage{geometry,cite}                
\usepackage{graphicx}
\usepackage{amssymb,amsmath}
\usepackage{graphicx}

\newcommand{\qbinom}[2] {{#1 \brack #2}_q}

\usepackage{epstopdf}
\newtheorem{thm}{Theorem}[section]
\newtheorem{cor}[thm]{Corollary}
\newtheorem{lem}[thm]{Lemma}
\newtheorem{prop}[thm]{Proposition}
\theoremstyle{definition}
\newtheorem{defn}[thm]{Definition}
\theoremstyle{remark}
\newtheorem{rem}[thm]{Remark}
\newtheorem{exa}[thm]{Example}
\numberwithin{equation}{section}

\newcommand{\be}{\begin{equation}}
\newcommand{\ee}{\end{equation}}

\newcommand{\bea}{\begin{eqnarray}}
\newcommand{\eea}{\end{eqnarray}}

\begin{document}

\title[]{ Odd and even  $q$-type Lidstone polynomial sequences }%
\author{Zeinab S.I. Mansour and Maryam  AL-Towailb}%
\address{Z. Mansour, Department of Mathematics, Faculty of Science, Cairo University, Giza, Egypt.}%
\email{zeinab@sci.cu.edu.eg}%
\address{M. AL-Towailb, Department of Computer Science and Engineering, King Saud University, Riyadh, KSA}%
\email{mtowaileb@ksu.edu.sa}%


\begin{abstract}
In this paper, we introduce two types of general classes of even and odd $q$-Lidstone polynomial sequences.  We prove essential properties related to them like the matrix and determinate form representation, the generating function, recurrence relation, and conjugate sequences. Some illustrative examples are included.
\end{abstract}



\maketitle

\section{\bf{Introduction}}

In 1929, Lidstone \cite{Lidstone} introduced a generalization of Taylor' series
that approximates a certain function in a neighborhood of two points insted of one. That is,
\[f(z)=\sum_{n=0}^{\infty}\Big[f^{(2n)}(1)\Lambda_n(z)-f^{{(2n)}}(0)\Lambda_{n}(z-1)\Big],\]
where $\Lambda_n$  is the Lidstone polynomials of degree $2n+1$ that given by
\begin{equation*}
\Lambda_n(z)=\frac{2^{2n+1}}{(2n+1)!}B_{2n+1}(\frac{1+z}{2}) ,
\end{equation*}
and $B_{n}(z)$ is the Bernoulli polynomials of degree $n$ (see \cite{Whittaker}).
 Recently, the theory of Lidstone polynomials has become an important investigation area where some results and number of problems was discussed in \cite{Agarwal1,Agarwal2,Agarwal3,Agarwal4,Agarwal5,Costabile1,Costabile2,Costabile3,Costabile4}. Ismail and Mansour constructed a $q$-Lidstone expansion theorem in the form
\[f(z)=\sum_{k=0}^{\infty}\Big[D_{q^{-1}}^{2n}f(1)B_n(z)-D_{q^{-1}}^{2n}f(0)A_n(z)\Big]. \]
Here $f$ is an entire function satisfying certain growth condition and the polynomials $A_n$ and $B_n$ defined by
\begin{equation}\label{A_n+B_n}
A_n(z) = \frac{2^{2n+1}}{[2n+1]_q!}B_{2n+1}(z/2;q), \;
    B_n(z) = \large{\eta}^1_{q^{-1}} A_n(z), 
\end{equation}
where $\large{\eta}^y_{q^{-1}}$ denotes the $q$-translation operator defined by
\begin{equation*}{\large{\eta}_{q^{-1}}^y} z^n
=q^{\frac{n(n-1)}{2}}z^n(-y/z  ;q^{-1})_n=y^n(-z/y  ;q)_n,\end{equation*}
and $B_{n}(z;q)$ is the $q$-Bernoulli polynomials defined by the
generating function
\begin{equation}\label{B(z;q)}
\dfrac{t\,E_q(zt)}{E_q(t/2)e_q(t/2)-1}=\sum_{n=0}^{\infty}B_n(z;q)
\frac{t^n}{[n]_q!},
\end{equation}
where $E_q(z)$ and $e_q(z)$ are the $q$-exponential functions defined by Jackson, cf. e.g. \cite{GR,Jacson 1},
\begin{equation}\label{E(z)&e(z)} E_q(z):= \sum_{j=0}^\infty q^{j(j-1)/2}\,\frac{z^j}{[j]_q!}; \, z\in \mathbb{C} \, \mbox{ and } \,
  e_q(z):= \sum_{j=0}^\infty \frac{z^j}{[j]_q!}; \, |z|< 1.\end{equation}
 Analogously,  Al-Towailb \cite{AL-Towaileb} introduced  the  $q$-Lidstone expansion:
\[f(z)=\sum_{n=0}^{\infty}\left[D_{q^{-1}}^{2n}f(1)\,M_n(z)- D_{q^{-1}}^{2n+1}f(0)\, N_{n}(z)
\right], \]
where  $M_n(z)$ and $N_n(z)$ are $q$-Lidstone polynomials defined by
\begin{equation}\label{M_n+N_n}
\left\{
  \begin{array}{ll}
M_n(z):= \frac{2^{2n}}{[2n]_q!}\large{\varepsilon}_{q^{-1}}^1E_{2n}(z/2;q),  & \hbox{} \\\\
    N_{n}(z)= \frac{2^{2n+2}}{[2n+1]_q!}E_{2n+1}(z/2;q). & \hbox{}
  \end{array}
\right.
\end{equation}
Here, $\large\varepsilon_{q^{-1}}^y$ a $q$-translation operator defined by
\begin{equation*}
{\large\varepsilon_{q^{-1}}^y} x^n
=x^n(-y/x  ;q^{-1})_n=q^{-\frac{n(n-1)}{2}}y^n(-x/y  ;q)_n
= q^{-\frac{n(n-1)}{2}}\large{\varepsilon}_{q}^x y^n,\end{equation*}
and $E_{n}(z;q)$ the $q$-Euler polynomials introduced by  Ismail and Mansour \cite{Ismail and Mansour}:
\begin{equation}\label{Euler Second kind}
\dfrac{2\,E_q(zt)}{e_q(t/2)E_q(t/2)+1}=\sum_{n=0}^{\infty}E_n(z;q)
\frac{t^n}{[n]_q!}.
\end{equation}

\noindent In \cite{Mansour and AL-Towaileb 2} the authors construct the complementary $q$-Lidstone  polynomials $\nu_n(z)$ and $\tau_n(z)$ of degree $2n$  satisfying 
\begin{equation}\label{nu and tau}
\left\{
  \begin{array}{ll}
    \nu_0(z)=1=\tau_0(z), & \hbox{} \\
    D_{q^{-1}} \nu_n(0)= D_{q^{-1}} \tau_n(1)=0, & \hbox{} \\
    D^2_{q^{-1}} \tau_n(z)= \tau_{n-1}(z) \, \mbox{ and } D^2_{q^{-1}} \nu_n(z)= \nu_{n-1}(z). & \hbox{}
  \end{array}
\right.
\end{equation}

\noindent Some applications related to $q$-Lidstone expansion theorem  are studied in \cite{Mansour and AL-Towaileb, Mansour and AL-Towaileb 2}.
\vskip3mm
In this paper, we will present general classes of $q$-polynomial sequences including \eqref{A_n+B_n}, \eqref{M_n+N_n} and \eqref{nu and tau}. For this,
we consider a sequence of $q$-polynomials
\begin{equation*}
L_0(z;q), L_1(z;q),\ldots, L_n(z;q), L_{n+1}(z;q), \ldots \quad (q\neq 0, \,  n\in \mathbb{N})
\end{equation*}
such that $L_n(z;q)$ satisfies the $q$-difference equation
\begin{equation}\label{qLPS1}
D^2_{q} L_n(z;q)= a_n\, L_{n-1}(z;q) \quad \mbox{ or } \quad  D^2_{q^{-1}} L_n(z;q)= a_n\, L_{n-1}(z;q),
\end{equation}
$(a_n\in\mathbb{R})$. We call $(L_n)_n$ a $q$-Lidstone Lidstone-type polynomial sequence ($q$LPS).

\noindent Note that when $q\rightarrow 1$, Equation \eqref{qLPS1} reduces to
\begin{equation*}
\frac{d^2}{dz^2} L_n(z)= a_n\, L_{n-1}(z),
\end{equation*}
so that we may think of this sequence as a generalization of Lidstone-type polynomial sequence (see \cite{Costabile,Costabile-Gualtieri-Napoli}).

\noindent Our aim is to study  some classes of $q$-Lidstone sets $\{L_n(z;q)\}_n$ when $q\neq 1$. More precisely, We will consider the  $q$-Lidstone polynomials of odd and even degrees in  the two types and introduce an algebraic approach to these polynomials. We follow the notations and terminologies in \cite{AMbook,GR}.
\vskip3mm

\noindent This article  is organized as follows: In the next section, we define the class of odd $q$-Lidstone polynomial sequences of type I, and we derive essential formulas  including matrix and determinate forms, the generating function, recurrence relation, and conjugate sequences. Then, in Section 3, we consider and study the class of even $q$-Lidstone polynomial sequences of type I. In Section 4, an analogy with the first class of  $q$-Lidstone polynomial sequences, we introduce $q$-Lidstone polynomial sequences of type II with theoretical results and some properties related to them. Finally, in  section 5, we give some illustrative examples.

\section{\bf{Odd $q$-Lidstone Polynomial Sequences of Type I}}\label{Sec:odd}

Throughout this paper, unless otherwise is stated,  $q$ is a positive number less than 1.  In this section, we shall define and study the  odd $q$-Lidstone polynomial sequences of Type I, which are in general satisfy the $q$-difference equation:   

\begin{equation}\label{qLPS}
D^2_{q} L_n(z;q)= a_n\, L_{n-1}(z;q), \quad a_n\in\mathbb{R} \, (n\in \mathbb{N}).
\end{equation}
\begin{defn}\label{qOLS}
 An  odd $q$-Lidstone polynomial sequence of type I ($q$OLS-I) is  a  set of polynomial sequences satisfying the second-order  $q$-difference equation
\begin{equation}\label{Def. of qOLS}
\left\{
  \begin{array}{ll}
D_{q}^{2}p_n(z;q)= [2n+1]_q\,[2n]_q\, p_{n-1}(z;q)& (n\in\mathbb{N}), \\\\
 p_n(0)=0 \,\; (n\in \mathbb{N}_0), \quad  p_0(z;q)=\alpha_0 z \; \alpha_0\in \mathbb{R}\backslash \{0\}. & \hbox{}
  \end{array}
\right.
\end{equation}
\end{defn}
One can verify that  $p_n(z;q)$ is a polynomial of degree $2n+1$ for all $n\in \mathbb{N}_0$.
\vskip2mm
\noindent A  series representation of the set $q$OLS-I is given in the following proposition.
\begin{prop}
A  sequence of polynomials  $\{p_n(z;q)\}_n$ is in the class  ( $q$OLS-I)  if and only if there exists a sequence $\{ \alpha_{2k}\}_{k\geq 0}$ of real numbers such that $\alpha_0 \neq 0$,  and
\begin{equation}\label{formala of qOLS}
\begin{split}
p_n(z;q)=& \sum_{k=0}^n \qbinom{2n+1}{2k+1}\, \frac{\alpha_{2k}}{[2(n-k)+1]_q}\, z^{2(n-k)+1}\\
=&\sum_{k=0}^n \qbinom{2n+1}{2k+1}\, \frac{\alpha_{2(n-k)}}{[2(n-k)+1]_q}\, z^{2k+1}.
\end{split}
\end{equation}
\end{prop}
\begin{proof}
Let $\{p_n(z;q)\}_n \in q$OLS-I. Then, there exists a constant $\alpha_0\neq 0$ such that $p_0(z;q)=\alpha_0 z$. Therefore,
\begin{equation*}
(D^2_{q} p_1)(z;q)= c_1 z, \quad c_1= [2][3]\,\alpha_0,
\end{equation*}
and then
\begin{equation*}
D_{q} p_1(z;q)= c_2 z^2+\alpha_2, \, c_2, \alpha_2 \in \mathbb{R}.
\end{equation*}
 This implies $D_{q} p_1(0)= \alpha_2$ and by induction, we can set \begin{equation}\label{D_qp_n}
D_{q} p_n(0;q)= \alpha_{2n}, \quad \alpha_{2n} \in \mathbb{R}.
\end{equation}
Now, assume that $p_n(z;q)=\sum_{k=0}^n \alpha_k^{(n)} z^{2(n-k)+1}$. Then,
\begin{equation}\label{Deff. Eq.1 }
D_{q}^{2}p_n(z;q)= \sum_{k=0}^{n-1} \alpha_k^{(n)} [2(n-k)+1]_q[2(n-k)]\, z^{2(n-k)-1}.
\end{equation}
According to the $q$-difference equation in \eqref{Def. of qOLS} we have
\begin{equation}\label{Deff. Eq.2 }
D_{q}^{2}p_n(z;q)=  [2n+1]_q\,[2n]_q\sum_{k=0}^{n-1} \alpha_k^{(n-1)}\, z^{2(n-k)-1}.
\end{equation}
From \eqref{Deff. Eq.1 } and \eqref{Deff. Eq.2 }, we get
\begin{equation*}
\prod_{n=k+1}^m \frac{\alpha_k^{(n)}}{\alpha_k^{(n-1)}}= \prod_{n=k+1}^m \frac{[2n+1]_q[2n]_q}{[2(n-k)+1]_q[2(n-k)]}.
\end{equation*}
Consequently, we obtain
\begin{equation}
\alpha_k^{(m)}= \frac{[2m+1]!}{[2(m-k)+1]![2k+1]_q!}\, \alpha_k^{(k)}= \qbinom{2m+1}{2k+1}\, \frac{\alpha_k^{(k)}}{[2(m-k)+1]},
\end{equation}
where $\alpha_k^{(k)}$ is the coefficient of $z$ in $p_k(z;q)$. Using \eqref{D_qp_n}, we can replace  $\alpha_k^{(k)}$ by $\alpha_{2k}$ and then we get the result in \eqref{formala of qOLS}.
 On the other hand, if \eqref{formala of qOLS} is satisfied and $\alpha_0\neq 0$, we get easily \eqref{Def. of qOLS} which complete the proof.
\end{proof}

\vskip .5 cm 

\begin{rem}
From \eqref{formala of qOLS}, we obtain
\begin{description}
  \item[i] $p_n(z;q)$ an odd function for all $n\in \mathbb{N}$;\\
  \item[ii] $\{z^{2n+1}\}_n\in qOLS$ and  $qOLS \subset \tilde{P}$, where $\tilde{P}=$ span$\{z^{2j+1}| \, j\in \mathbb{N}\}$;\\
\item[iii] $ \int_0^1 p_n(z;q)\, d_qz= [2n+1]_q! \sum_{k=0}^n  \frac{\alpha_{2(n-k)}}{[2k+2]!\, [2(n-k)+1]_q!}, \, n\in\mathbb{N}$.
\end{description}
\end{rem}

\vskip .5 cm 

\begin{prop}\label{Prop.1}
Let $n\in \mathbb{N}$ and $\{p_n(z;q)\}_n\in q$OLS-I. Then
\begin{enumerate}
         \item $D_{q}^{2m} p_n(z;q)= \frac{[2n+1]_q!}{[2(n-m)+1]!}\, p_{n-m}(z;q), \quad  m=1,2,\ldots,n; $\\
         \item $D_{q}^{2m+1} p_n(z;q)= \frac{[2n+1]_q!}{[2(n-m)+1]!}\, D_{q}p_{n-m}(z;q), \quad  m=1,2,\ldots,n;$\\
         \item $D_{q}^{2m} p_n(0)=0; \, D_{q}^{2m+1} p_n(0)= \frac{[2n+1]_q!}{[2(n-m)+1]!}\, \alpha_{2(n-m)}, \quad  m=1,2,\ldots,n$.
         \end{enumerate}
\end{prop}
\begin{proof}
The proof follows immediately from \eqref{Def. of qOLS} and \eqref{formala of qOLS} by induction.
\end{proof}
\vskip3mm
\subsection{The matrix form } ~\par
\vskip1mm

In the following, we present the matrix form $A_q$ of the odd $q$-Lidstone polynomials of type I and then, we analyze the structure of this matrix.
\vskip3mm
\noindent Recall that a matrix $M=[m_{ij}] \, (i,j\in \mathbb{N}_0)$ is infinite lower triangular if $m_{ij}=0$ whenever $j>i$. Denote by $\mathcal{L}$ the set of all lower triangular matrices.
\vskip2mm
\noindent A matrix $T=[a_{ij}] \in \mathcal{L}$ is a Toeplitz if and only if $a_{ij}=a_{i+mj+m}$ for all nonnegative integers $i,j,m$. If $T$ is
a Toeplitz matrix, define the sequence $a_n:=a_{n0}$ for $n\geq 0$,  i.e.  $a_{ij}=a_{i-j}$ for $i\geq j$. That is,
\begin{equation*}
T= \left[
     \begin{array}{cccccccccc}
        & \vdots & \vdots & \vdots & \vdots \\
       \ldots & a_0 & 0 & 0 & 0 & \ldots \\
       \ldots & a_1 & a_0 & 0 & 0 & \ldots \\
       \ldots & a_2 & a_1 & a_0 & 0 & \ldots \\
      \ldots & a_3 & a_2 & a_1 & a_0 & \ldots \\
         & \vdots & \vdots & \vdots & \vdots \\
     \end{array}
   \right].
\end{equation*}

Verde-Star \cite{Verde}  introduced Lemma \ref{the prodect of two matrices} and Lemma \ref{invers of T} below.  

\begin{lem}\label{the prodect of two matrices}
Let $A, B \in \mathcal{L}$. Then the product $C=AB$ is well defined of $\mathcal{L}$ and
\begin{equation*}
c_{ik}= \sum_{j=k}^i a_{ij}\, b_{jk}, \, i\geq k.
\end{equation*}
\end{lem}
\begin{lem}\label{invers of T}
Let $T_a$ be a Toeplitz matrix in $\mathcal{L}$ defined by
\begin{equation*}
T_a= [a_{ij}]:= \left\{
                  \begin{array}{ll}
                    a_{i-j},  & \hbox{$i\geq j$;} \\
                    0, & \hbox{$i< j$.}
                  \end{array}
                \right.
\end{equation*}
If $(T_a)^{-1}= T_b:= [b_n] \, (n=i-j)$, then
\begin{equation*}
b_n= \frac{(-1)^n}{a_0^{n+1}}\, \mbox{det} \left[
                                             \begin{array}{ccccc}
                                               a_1 & a_0 &  &  &  \\
                                               a_2 & a_1 & a_0 &  &  \\
                                               \vdots & \vdots & \vdots & \ddots &  \\
                                               a_{n-1} & a_{n-2} & a_{n-3} & \ldots & a_0 \\
                                               a_n & a_{n-1} & a_{n-2} & \ldots & a_1 \\
                                             \end{array}
                                           \right].
\end{equation*}
\end{lem}

\begin{defn}\label{matrix form}
The  odd $q$-Lidstone-type I matrix is an infinite lower triangular matrix $A_q=[a_{ij}]$ ($i,j\in \mathbb{N}_0$) with
\begin{equation}\label{elements of A}
a_{ij}= \qbinom{2i+1}{2j+1}\, \frac{\alpha_{2(i-j)}}{[2(i-j)+1]_q}, \, i\geq j,
\end{equation}
where $\{\alpha_{2k}\}_{k\geq 0}$ a sequence of real numbers and $\alpha_0\neq 0$.
\end{defn}
\vskip .2 cm

\begin{rem}
According to Definition \ref{matrix form}, Formula \eqref{formala of qOLS} can be written in the  matrix form
\begin{equation}\label{matrix of olp}
P_q= A_qZ_q,
\end{equation}
where $P_q$ and $Z_q$ are two vectors defined by
\begin{equation*}
P_q= [p_0(z;q), p_1(z;q),\ldots, p_n(z;q), \ldots]^T, \quad Z_q=[z, z^3, \ldots, z^{2n+1}, \ldots]^T.
\end{equation*}
Moreover, if we set $A_{q,n}=[a_{ij}]$ such that $j=0,1,...,i$, $i=0,1, ...,n$ for $n\in \mathbb{N}$ and $a_{ij}$ defined in \eqref{elements of A}, then we have a sequence $(A_{q,n})_{n}$ of the principle submatrices of order $n$ of $A_q$ which satisfy
\begin{equation}\label{finite matrix of olp}
P_{q,n}= A_{q,n}Z_{q,n},
\end{equation}
where
\begin{equation}\label{vector spaces P and Z}
P_{q,n}= [p_0(z;q), p_1(z;q),\ldots, p_n(z;q)]^T, \quad Z_{q,n}=[z, z^3, \ldots, z^{2n+1}]^T.
\end{equation}
\end{rem}

\vskip .2 cm 

\begin{rem}
From \eqref{vector spaces P and Z}, we get
\begin{equation*}
D_{q}^2P_{q,n}= [D_{q}^2p_0, D_{q}^2p_1,\ldots, D_{q}^2p_n]^T.
\end{equation*}
So, if we denote by $\mathcal{C}=[c_{ij}] \, (i,j=0,1,2,\ldots,n)$ to the derivation matrix for $P_{q,n}$, that is
\begin{equation*}
D_{q}^2P_{q,n}= \mathcal{C}P_{q,n}, \quad n\in \mathbb{N}_0,
\end{equation*}
then According to \eqref{Def. of qOLS}, we obtain
\begin{equation*}
c_{ij}=\left\{
    \begin{array}{ll}
      [2i+1]_q[2i]_q , & \hbox{$i = j+1$;}  \\
0, & \hbox{otherwise}.
    \end{array}
  \right.
\end{equation*}
\end{rem}

\vskip3mm

\noindent We define the Toeplitz matrix $T_{q\alpha}$ as the matrix in $\mathcal{L}$ whose $(i,j)$ entry is
\begin{equation*}
t_{ij}^{\alpha}=\frac{\alpha_{2(i-j)}}{[2(i-j)+1]_q!}
\end{equation*}
for $ i\geq j$, and zero otherwise. Also, we denote  $\mathcal{D}$ to the diagonal matrix with entries $d_{ii}= [2i+1]_q!$.
\begin{prop}
The odd $q$-Lidstone-type I matrix $A_q$ can be factorized as
\begin{equation}\label{factorize A_q}
 A_q= \mathcal{D}T_{q\alpha} \mathcal{D}^{-1}.
\end{equation}
\end{prop}
\begin{proof}
According to Lemma \ref{the prodect of two matrices}, the product $\mathcal{D}T_{q\alpha} \mathcal{D}^{-1}$ is well-defined and we easily get the result.
\end{proof}
\begin{prop}\label{(A_q)^{-1}}
The  odd $q$-Lidstone-type I matrix $A_q$ is invertible and
\begin{equation}\label{2factorize A_q}
 (A_q)^{-1}= \mathcal{D}T_{q\beta} D^{-1},
\end{equation}
where $(\beta_{2n})_n$ is the numerical sequence satisfying 
\begin{equation}\label{beta sequence}
 \sum_{j=0}^n \frac{\beta_{2j}\alpha_{2(n-j)}}{[2j+1]_q![2(n-j)+1]_q!}= \delta_{n0}, \quad n\in \mathbb{N}_0,
\end{equation}
with $\delta_{nj}$ the Kronecker's delta.
 \end{prop}
\begin{proof}
It follows directly by using \eqref{factorize A_q} and calculating $(T_{q\alpha})^{-1}$ from the result in Lemma \ref{invers of T}.
\end{proof}
\begin{rem}
Equation \eqref{beta sequence} can be considered as an infinite linear system which determines the numerical sequence $(\beta_{2n})_n$ and,
according to Cramer's rule, the first $m + 1$ equations give
\begin{equation}\label{matrix of beta}
\begin{split}
\beta_0=&\frac{1}{\alpha_0},\\
\beta_{2n}=& (-1)^n\frac{[3]![5]! \ldots [2n+1]_q!}{\alpha_0^{n+1}}\\
\times& \mbox{det}\left[
          \begin{array}{ccccc}
             \frac{\alpha_2}{[3]!} & \frac{\alpha_0}{[3]_q!} & 0 &  \ldots & 0 \\
             \frac{\alpha_4}{[5]!} & \frac{\alpha_2}{[3]_q![3]_q!} & \frac{\alpha_0}{[5]_q!} & 0  & 0 \\
            \vdots & \vdots & \vdots & \ddots & \vdots \\
            \vdots & \vdots & \vdots & \ddots  & \vdots \\
             \frac{\alpha_{2(n-1)}}{[2n-1]_q!} & \frac{\alpha_{2(n-2)}}{[2n-3]_q![3]_q!} & \frac{\alpha_{2(n-3)}}{[2n-5]_q![5]_q!} & \ldots & \frac{\alpha_{0}}{[2n-1]_q!} \\
           \frac{\alpha_{2n}}{[2n+1]_q!}  & \frac{\alpha_{2(n-1)}}{[2n-1]_q![3]_q!} & \frac{\alpha_{2(n-2)}}{[2n-3]_q![5]_q!} & \ldots & \frac{\alpha_{2}}{[2n-1]_q![3]_q!} \\
          \end{array}
        \right].
\end{split}
\end{equation}
\end{rem}
\vskip5mm

Now, we consider the polynomials
\begin{equation}\label{formala of qOLS of beta}
\hat{p}_n(z;q)= \sum_{k=0}^n \qbinom{2n+1}{2k+1}\, \frac{\beta_{2(n-k)}}{[2(n-k)+1]_q}\, z^{2k+1},
\end{equation}
where $(\beta_{2n})_n$ is  the numerical sequence defined in \eqref{matrix of beta}. Note that $\{\hat{p}_n(z;q)\}_n\in q$OLS-I.

\begin{defn}
The two sequences
$\{p_n(z;q)\}_n$ and $\left(\hat{p}_n(z;q)\right)_n$ defined in  \eqref{formala of qOLS} and  \eqref{formala of qOLS of beta}, respectively, are called conjugate odd $q$-Lidstone-type sequences.
\end{defn}

\vskip5mm
\noindent We denote $B_q=[b_{ij}], \, i,j\in \mathbb{N}_0$ the infinite lower triangular matrix with
\begin{equation*}
b_{ij}=  \qbinom{2i+1}{2j+1}\, \frac{\beta_{2(i-j)}}{[2(i-j)+1]_q}, \quad i\geq j.
\end{equation*}
Set $\widehat{P}_q= [\hat{p}_0(z;q), \hat{p}_1(z;q),\ldots, \hat{p}_n(z;q), \ldots]^T$. Then, we have the  matrix forms
\begin{equation}\label{matrix of olp-beta}
\widehat{P}_q= B_qZ_q,
\end{equation}
and for $n\in \mathbb{N}$,
\begin{equation}\label{matrix of olp-beta-n}
\widehat{P}_{q,n}= B_{q,n}Z_{q,n}.
\end{equation}

\begin{prop}\label{conjugate $q$-odd Lidstone-type sequences}
The  sequences $\{p_n(z;q)\}_n$ and $\left(\hat{p}_n(z;q)\right)_n$  are  conjugate odd $q$-Lidstone-type sequences if and only if
\begin{equation*}
P_q= A_q^2\widehat{P}_q, \quad  \widehat{P}_q= B_q^2P_q,
\end{equation*}
and for $n\in \mathbb{N}$,
\begin{equation*}
P_{q,n}= A_{q,n}^2\widehat{P}_{q,n}, \quad  \widehat{P}_{q,n}= B_{q,n}^2P_{q,n}.
\end{equation*}
\end{prop}
\begin{proof}
The proof  follows directly from  \eqref{matrix of olp}, \eqref{finite matrix of olp}, \eqref{matrix of olp-beta} and \eqref{matrix of olp-beta-n}, and taking into account $A_q^{-1}=B_q$.
\end{proof}
\begin{rem}
From Proposition \ref{conjugate $q$-odd Lidstone-type sequences}, we can write
\begin{equation*}
p_n(z;q)= \sum_{k=0}^n \tilde{a}_{nk}\, \hat{p}_k(z;q), \quad
\hat{p}_n(z;q)= \sum_{k=0}^n \tilde{b}_{nk}\, p_k(z;q) \, (n\in \mathbb{N}_0),
\end{equation*}
where $\tilde{a}_{nk}$ and $\tilde{b}_{nk}$, $k=0,\ldots, n$, are the elements of the matrices $A_{q,n}^2$ and $B_{q,n}^2$, respectively.
\end{rem}

\vskip3mm
\subsection{Recurrence relations and $q$-difference  equations}~\par
\vskip1mm

\noindent  We start this subsection by deriving some recurrence relations for the sequence of odd $q$-Lidstone polynomials of type I.

\begin{thm}\label{First recu. relation Th.}
Let $\{p_n(z;q)\}_n\in q$OLS-I. Then
\begin{equation}\label{First recu. relation of qOLS}
p_n(z;q)= \frac{1}{\beta_0}\Big[z^{2n+1}-\sum_{k=0}^{n-1} \qbinom{2n+1}{2k+1}\, \frac{\beta_{2(n-k)}}{[2(n-k)+1]_q}\, p_k(z;q)\Big], \, n\in \mathbb{N}_0,
\end{equation}
where $(\beta_{2n})_n$ is the numerical sequence defined as in \eqref{matrix of beta}.
\end{thm}

\begin{proof}
Assume that $A_q$ is  an odd $q$-Lidstone-type I matrix. Then, from \eqref{finite matrix of olp} we have
\begin{equation}\label{Z=BP}
Z_{q,n}= B_{q,n}P_{q,n},
\end{equation}
where $Z_{q,n}$ and $P_{q,n}$ are defined as in \eqref{vector spaces P and Z} and $B_{q,n}=A_{q,n}^{-1}$.
Therefore, for $i=0,1,...,n$ we obtain
\begin{equation}\label{First recu. relation of qOLS-1}
z^{2n+1}= \sum_{k=0}^{n} \qbinom{2n+1}{2k+1}\, \frac{\beta_{2(n-k)}}{[2(n-k)+1]_q}\, p_k(z;q),
\end{equation}
and then we get the result.
\end{proof}
\vskip3mm

Note that Relation \eqref{First recu. relation of qOLS-1} can be considered as infinite linear system in the unknowns $(p_n(z;q))_n$.  In the following theorem, we use 
Cramer's rule to solve the first $n+1$ equations of this system, and then, we obtain a first determinate form of odd $q$-Lidstone polynomial sequence.

\begin{thm}
If  $\left(p_n(z;q)\right)_n$ is in the class  $q$OLS-I, then
\begin{equation}\label{a first determinate form}
\begin{split}
p_0(z;q)=&\frac{1}{\beta_0}z,\\
p_{n}(z;q)=& \frac{(-1)^n}{[3]_q![5]_q! \ldots [2n-1]_q!\,\beta_0^{n+1}}\\
\times& \mbox{det}\left[
          \begin{array}{cccccc}
             z & z^3 & z^5 & \ldots & z^{2n-1} & z^{2n+1} \\
             \beta_0 & \beta_2 & \beta_4 & \ldots & \beta_{2(n-1)} & \beta_{2n} \\
            0 & [3]!\beta_0 & \frac{[5]_q!}{[3]_q!}\beta_2 & \ldots & \frac{[2n-1]_q!}{[2n-3]_q!}\beta_{2(n-2)}& \frac{[2n+1]_q!}{[2n+3]!}\beta_{2(n-1)} \\
            \vdots & \ddots & \ddots & & \vdots & \vdots \\
\vdots & \ddots & \ddots & & \vdots & \vdots \\
 0 & \ldots & \ldots & & [2n-1]_q!\beta_0 & \frac{[2n+1]_q!}{[3]!}\beta_2
          \end{array}
        \right],
\end{split}
\end{equation}
where $(\beta_{2n})_n$  is the numerical sequence defined as in \eqref{matrix of beta}.
\end{thm}
\vskip2mm

The following  theorem  gives    a recurrence relation for the conjugate sequence $\{\hat{p}_n(z;q)\}_n$.

\begin{thm}\label{Thm:z}
Let $\{p_n(z;q)\}_n$ be an odd $q$-Lidstone polynomial sequence of type I. Then, the conjugate odd $q$-Lidstone polynomial sequence $\{\hat{p}_n(z;q)\}_n$,  satisfies the recursive relation
\begin{equation}\label{Rec.Rel.conj.}
\hat{p}_n(z;q)= \frac{1}{\alpha_0}\Big[z^{2n+1}-\sum_{k=0}^{n-1} \qbinom{2n+1}{2k+1}\, \frac{\alpha_{2(n-k)}}{[2(n-k)+1]_q}\, \hat{p}_k(z;q)\Big],
\end{equation}
where  $(\alpha_{2n})_{n\geq 0}$  is the numerical sequence defined as in \eqref{formala of qOLS}.
Moreover,  $\left(\hat{p}_n(z;q)\right)_n$ can be expressed in a determinate form similar to \eqref{a first determinate form} with $\alpha_{2k}$ instead of
$\beta_{2k}$, for $k=0, 1, ..., n$ and $n\in \mathbb{N}_0$.
\end{thm}

\begin{proof}
The proof is similar to the proof of Theorem \ref{First recu. relation Th.} and is omitted.
\end{proof}
\vskip3mm

Now, we determine another recurrence relation by using the production matrix. For this, we recall the definition and some properties of this matrix (for more details, see \cite{Deutsch}).
 \begin{defn}
Let $A$ be a nonsingular  infinite lower triangular matrix. The production matrix $\Pi_A$ of $A$ is defined by
\begin{equation*}
\Pi_A= A^{-1}\bar{A},
\end{equation*}
where $\bar{A}$ is the matrix $A$ with its first row removed.
\end{defn}
\begin{prop}
Let $A$ be an infinite lower triangular matrix and $B$ be the inverse matrix. Then, the production matrices  $\Pi_A$ and
 $\Pi_B$ of $A$ and $B$, respectively, satisfy
\begin{equation}\label{prod.}
\Pi_B A= AD \quad \mbox{ and } \quad \Pi_A B= BD,
\end{equation}
where $D=[\delta_{(i+1)j}] \, (i,j\in \mathbb{N}_0)$ and $\delta_{ij}$ is the Kronecker's delta.
\end{prop}

\begin{lem}
Let $A_q=[a_{ij}] \, (i,j\in \mathbb{N}_0)$ be an odd $q$-Lidstone-type I matrix, $B_q=[b_{ij}]$ be the inverse of $A_q$,  and $\Pi_B=[\pi_{ij}]$ be the production matrix of $B_q$. Then
\begin{equation}\label{elements of prod.}
\begin{split}
\pi_{ij}=& \sum_{n=0}^i a_{in} b_{(n+1)j}\\
=&\left\{
    \begin{array}{ll}
      \alpha_0\beta_2, & \hbox{$i=j=0$,} \\
      0, & \hbox{$j > i+1$,}  \\
     \sum_{n=0}^{i-j+1} \qbinom{2i+1}{2(n+j)-1}\, \frac{\beta_{2n}\alpha_{2(i-j-n)+2}[2(n+j)+1]_q!}{(2(i-j-n)+3)[2j+1]_q![2n+1]_q!}, & \hbox{otherwise,}
    \end{array}
  \right.
\end{split}
\end{equation}
where $(\alpha_{2n})_n$ and $(\beta_{2n})_n$ the numerical sequences defined as in \eqref{beta sequence}.
\end{lem}
\begin{proof}
From \eqref{prod.}, we have $\Pi_B=ADB$. Thus, $\pi_{ij}= \sum_{n=0}^i a_{in} b_{(n+1)j}$, and by Proposition \ref{(A_q)^{-1}},  we get the result.
\end{proof}

\begin{thm}\label{Second recu. relation thm}
Let $\{p_n(z;q)\}_n$ be  in the class  $q$OLS-I. If $A_q$ is the odd $q$-Lidstone-type I matrix related to  $\{p_n(z;q)\}_n$, and
  $\Pi_q=[\pi_{ij}]\, (i,j\in \mathbb{N}_0)$ is the production matrix of $A^{-1}_q=B_q$, then
\begin{equation}\label{Second recu. relation}
\begin{split}
p_0(z;q)=&\frac{1}{\beta_0}z, \\
p_{n+1}(z;q)=& \frac{1}{\pi_{n(n+1)}}\Big[z^{2}p_{n}(z;q)-\sum_{k=0}^{n} \pi_{nk} p_k(z;q)\Big], \, n\in \mathbb{N}_0.
\end{split}
\end{equation}
\end{thm}
\begin{proof}
From \eqref{matrix of olp} and \eqref{prod.}, we have $\Pi_q P_q=A_q(DZ_q)$. Since $DZ_q=[z^3,z^5,...]^T= z^2Z_q$, we obtain
\begin{equation}\label{Prod.2}
\Pi_q P_q=z^2A_qZ_q=z^2P_q.
\end{equation}
Consider the first $(n+1)$ equations in \eqref{Prod.2}, we have $\sum_{k=0}^{n+1} \pi_{nk} p_k(z;q)= z^{2}p_{n}(z;q)$, which is nothing else but \eqref{Second recu. relation}. This completes the proof of the theorem.  
\end{proof}
\begin{thm}
Let $\{p_n(z;q)\}_n\in q$OLS-I. Then
 \begin{equation*}
\begin{split}
p_0(z;q)=&\frac{1}{\beta_0}z,\\
p_{n+1}(z;q)=& \frac{(-1)^{n+1}p_0(z;q)}{\pi_{01}\pi_{12} \ldots \pi_{n(n+1)}}\\
\times& \mbox{det}\left[
          \begin{array}{cccccc}
           \pi_{00}-z^2 & \pi_{01} & 0 & \ldots & \ldots & 0 \\
             \pi_{10} & \pi_{11}-z^2 & \pi_{12} & \ldots & \ldots & 0 \\
            \pi_{20} & \pi_{21} & \pi_{22}-z^2 & \ddots & \ddots & 0 \\
            \vdots & \vdots & \vdots &  \ddots & \ddots & \vdots\\
\vdots & \vdots & \vdots & \ddots & \ddots & \pi_{(n-1)n} \\
 \pi_{n0} & \pi_{n1} & \pi_{n2} & \ldots & \ldots & \pi_{nn}-z^2
          \end{array}
        \right],
\end{split}
\end{equation*}
where $\pi_{ij}$ are  defined as in \eqref{elements of prod.}.
\end{thm}
\begin{proof}
According to Theorem \ref{Second recu. relation thm}, we have the linear system \eqref{Prod.2} which can be expressed in a matrix form as
 \begin{equation*}
\left[
          \begin{array}{ccccc}
           \pi_{01} & 0 & 0 & 0 & \ldots  \\
            \pi_{11}-z^2 & \pi_{12} & 0 & 0 & \ldots \\
            \pi_{21} & \pi_{22}-z^2 & \pi_{23} & 0 & \ldots  \\
            \pi_{31} & \pi_{32} & \pi_{33}-z^2 & \pi_{34} & \ldots \\
\vdots & \vdots & \vdots & \vdots & \ddots
          \end{array}
        \right] \left[
                  \begin{array}{c}
                    p_1 \\
                    p_2 \\
                    p_3 \\
                    p_4 \\
                    \vdots \\
                  \end{array}
                \right]
= p_0\left[
             \begin{array}{c}
               z^2-\pi_{00} \\
               -\pi_{10} \\
               -\pi_{20} \\
               -\pi_{30} \\
                \vdots \\
             \end{array}
           \right].
\end{equation*}
By using Cramer's rule, we get the solution of the first $n+1$ equations and then, we obtain the result.
\end{proof}
\vskip5mm

\noindent We end this section by proving that the odd $q$-Lidstone polynomials of type I satisfy some $q$-difference equations.
\begin{thm}
If $p_n(z;q)\, (n\in \mathbb{N}_0)$ is an odd $q$-Lidstone polynomial of type I, then it satisfies the following linear $q$-difference equations of order $2n$:
\begin{equation}\label{First Diff. Eq.}
\sum_{k=0}^n \frac{\beta_{2k}}{[2k+1]_q!}\, D_{q}^{2k}u(z)- z^{2n+1}=0,
\end{equation}
\begin{equation}\label{Second Diff. Eq.}
\sum_{k=1}^{n}  \frac{[2k+1]_q!}{[2n+1]_q!}\, \pi_{nk} D_{q}^{2(n-k+1)}u(z)-
D_{q}^2(z^2u(z))  +[2n+3][2n+2]\pi_{n(n+1)} u(z)=0.
\end{equation}
\end{thm}
\begin{proof}
From Proposition \ref{Prop.1}, we have
\begin{equation}\label{Diff.}
D_{q}^{2k} p_n(z;q)= \frac{[2n+1]_q!}{[2(n-k)+1]_q!}\, p_{n-k}(z;q).
\end{equation}
Substituting \eqref{Diff.} into \eqref{First recu. relation of qOLS}, we obtain
\begin{equation*}
p_{n+1}(z;q)= \frac{1}{\beta_0}\Big[z^{2n+3}-\sum_{k=0}^{n} \frac{[2n+3][2n+2]}{[2k+3]!} \,\beta_{2k+2}D_{q}^{2k} p_n(z;q)\Big].
\end{equation*}
Therefore, \begin{equation*}
\begin{gathered}
D_{q}^{2} p_{n+1}(z;q)= \\ \frac{1}{\beta_0}\Big[[2n+3][2n+2]z^{2n+1}- \sum_{k=0}^{n} \frac{[2n+3][2n+2]}{[2k+3]!}\,\beta_{2k+2}D_{q}^{2k+2} p_n(z;q)\Big].
\end{gathered}
\end{equation*}
Since \[D_{q}^{2} p_{n+1}(z;q)= [2n+2][2n+3]p_{n}(z;q),\;\mbox{ and}\; D_{q}^{2n+2} p_{n}(z;q)=0,\] we obtain \eqref{First Diff. Eq.}. 
 On the other hand, from Theorem \ref{Second recu. relation thm},  we have
\begin{equation}\label{S.RE.E.}
\pi_{n(n+1)}p_{n+1}(z;q)= z^{2}p_{n}(z;q)-\sum_{k=0}^{n} \pi_{nk} p_{k}(z;q).
\end{equation}
Taking the second $q$-derivative of \eqref{S.RE.E.}, we obtain
\begin{equation}\label{S.RE.E.2}
[2n+3][2n+2]\pi_{n(n+1)}p_{n}(z;q)= D_q^2(z^{2}p_{n}(z;q))-\sum_{k=1}^{n}  [2k+1]_q[2k]_q \pi_{nk} p_{k-1}(z;q).
\end{equation}
According to Proposition \ref{Prop.1}, 
\begin{equation}\label{S.RE.E.3}
D_{q}^{2(n-k+1)} p_n(z;q)= \frac{[2n+1]_q!}{[2k-1]!}\, p_{k-1}(z;q).
\end{equation}
 Substituting \eqref{S.RE.E.3} into \eqref{S.RE.E.2} yields \eqref{Second Diff. Eq.} and completes the proof.
\end{proof}

\vskip3mm
\subsection{Generating function}~\par
\vskip1mm

\noindent Our aim here is to calculate the generating function of the odd $q$-Lidstone polynomial sequence of type I.
\vskip2mm
Recall the $q$-trigonometric functions
 \[\text{Sin}_q (z):=\frac{E_q(iz)-E_q(-iz)}{2i},  \quad \sin_q(z)
=\frac{e_q(iz)-e_q(-iz)}{2i},\]
 \[\text{Cos}_q (z):=\frac{E_q(iz)+E_q(-iz)}{2},  \quad \cos_q(z)
=\frac{e_q(iz)+e_q(-iz)}{2i},\]
where $E_q(z)$ and $e_q(z)$ are the $q$-exponential functions defined in \eqref{E(z)&e(z)}.

\noindent Here, the $q$-analog of the hyperbolic functions $\sinh z$ and $\cosh z$ are defined for $z\in \mathbb{C}$ by
\be\label{sh+cosh}
 \text{sinh}_q(z):= -i \sin_q (iz)= \sum_{n=0}^{\infty}\frac{z^{2n+1}}{[2n+1]_q!}, \quad
 \text{cosh}_q(z):=  \cos_q(iz) = \sum_{n=0}^{\infty}\frac{z^{2n}}{[2n]_q!}.\ee
Also, $\text{Sinh}_q (z)=\text{sinh}_{1/q}(z)$ and $\text{Cosh}_q (z)=\text{cosh}_{1/q}(z)$.
\vskip 5mm

 Let $\{p_n(z;q)\}_n$ be the odd $q$-Lidstone polynomial sequence related to the numerical sequence $(\alpha_{2n})_{n\geq 0}$, and consider the following power series
\begin{equation}\label{power series g}
g_q(t)=\sum_{n=0}^\infty \alpha_{2n} \frac{t^{2n}}{[2n+1]_q!}.
\end{equation}

\begin{lem}\label{Lem.power series 1/g}
Let $g_q(t)$ be the power series defined in \eqref{power series g}. Then $\frac{1}{g_q(t)}$ is a well-defined function and has the series representation  
\begin{equation}\label{power series 1/g}
\frac{1}{g_q(t)}=\sum_{n=0}^\infty \beta_{2n} \frac{t^{2n}}{[2n+1]_q!},
\end{equation}
where $(\beta_{2n})_n$  is the numerical sequence defined as in \eqref{matrix of beta}.
\end{lem}

\begin{proof}
Since $\alpha_0\neq 0$, $g_q(t)$ is invertible. To prove the theorem, we need to prove that 
\begin{equation*}
\left(\sum_{n=0}^\infty \alpha_{2n} \frac{t^{2n}}{[2n+1]_q!}\right)\;\left(\sum_{n=0}^\infty \beta_{2n} \frac{t^{2n}}{[2n+1]_q!}\right)=1.
\end{equation*}
 By using the Cauchy product for power series, we obtain 
\[\begin{gathered}\left(\sum_{n=0}^\infty \alpha_{2n} \frac{t^{2n}}{[2n+1]_q!}\right)\;\left(\sum_{n=0}^\infty \beta_{2n} \frac{t^{2n}}{[2n+1]_q!}\right)\\=\sum_{n=0}^{\infty}t^{2n}\sum_{k=0}^{n}\frac{\beta_{2k}}{[2k+1]_q!} \frac{\alpha_{2(n-k)}}{[2(n-k)+1]_q!}.\end{gathered}\]
Thus,  \eqref{power series 1/g}  follows from \eqref{beta sequence}.

\end{proof}

\begin{thm}\label{Sec2:gf}
Let $\{p_n(z;q)\}_n$ and $\left(\hat{p}_n(z;q)\right)_n$  be the conjugate odd $q$-Lidstone-type I sequences. Then
\begin{eqnarray}\label{g1}
 g_q(t)\frac{\sinh_q (zt)}{t}= \sum_{n=0}^\infty p_n(z;q) \frac{t^{2n}}{[2n+1]_q!}, \\ \label{g2}
\frac{1}{g_q(t)} \frac{\sinh_q (zt)}{t}= \sum_{n=0}^\infty \hat{p}_n(z;q) \frac{t^{2n}}{[2n+1]_q!},
\end{eqnarray}
where $g_q(t)$ the function defined as in \eqref{power series g}.
\end{thm}
\begin{proof}
From \eqref{First recu. relation of qOLS-1}, we have
\begin{equation}\label{=recu.}
 \sum_{k=0}^{n} \qbinom{2n+1}{2k+1}\, \frac{\beta_{2(n-k)}}{[2(n-k)+1]_q}\, p_k(z;q)=z^{2n+1}\quad (n\in \mathbb{N}_0).
\end{equation}
Multiplying both sides of \eqref{=recu.} by $\frac{t^{2n+1}}{[2n+1]_q!}$ and adding on $n$, we obtain
\begin{equation*}
\sum_{n=0}^{\infty}\Big( \sum_{k=0}^{n} \qbinom{2n+1}{2k+1}\, \frac{\beta_{2(n-k)}}{[2(n-k)+1]_q}\, p_k(z;q)\Big)\frac{t^{2n+1}}{[2n+1]_q!} =\sum_{n=0}^{\infty}\frac{(zt)^{2n+1}}{[2n+1]_q!}.
\end{equation*}
Therefore,
\begin{equation*}
t\sum_{n=0}^\infty \beta_{2n} \frac{t^{2n}}{[2n+1]_q!}\sum_{n=0}^\infty p_n(z;q) \frac{t^{2n}}{[2n+1]_q!}= sinh_q(zt).
\end{equation*}
Now, using Lemma \ref{Lem.power series 1/g} we obtain \eqref{g1}.

\noindent Similarly, from \eqref{Rec.Rel.conj.}  we can derive the generating function of $ \hat{p}_n(z;q)$ and get Equation \eqref{g2}.
\end{proof}
\begin{exa}
Ismail and Mansour~\cite[Eq.(3.37)]{Ismail and Mansour} introduced the identity 
\begin{equation}\label{IM}
\sum_{k=0}^{n}(-1)^k 2^{2k}\frac{\beta_{2k}(q)}{[2k]_q!} \frac{T_{2n-2k+1}(q)}{[2n-2k+1]_q!}=\delta_{n,0},
\end{equation}
where $\delta_{n,0}$ is the Kroncker's delta.  It is worth noting that there was a small typo in \cite[Eq.(3.37)]{Ismail and Mansour},  which we have corrected in \eqref{IM}. 
  We  can take 
\[\beta_{2j}=(-1)^j 2^{2j}[2j+1]_q\beta_{2j}, \quad \alpha_{2j}=T_{2j+1}.\]
Then the polynomials $p_n(z;q))_n$ and $(\hat{p}_n(z;q))_n$ defined by 
\[p_n(z;q):=\sum_{k=0}^{n}\qbinom{2n+1}{2k+1} (-1)^k 2^{2k}\beta_{2k} z^{2n-2k+1},\]
\[\hat{p}_n(z;q):=\sum_{k=0}^{n}\qbinom{2n}{2k} \frac{T_{2n-2k+1}}{[2n-2k+1]_q} z^{2k+1},\]
are conjugate odd  $q$-Lidstone polynomials of type I. 
Since 
\[\text{Tan}_q t=\tan_qt=\sum_{n=0}^{\infty}T_{2n+1}\frac{t^{2n+1}}{[2n+1]_q!} ,\]
 \[t \text{Cot}_q t=t\cot_qt=\sum_{n=0}^{\infty}(-1)^n\beta_{2n}(q)\frac{ t^{2n}}{[2n]_q!},\]
 we can prove that 
 \[\frac{\tan_qt}{t}\frac{\sinh_q(zt)}{t}=\sum_{n=0}^{\infty}\hat{p}_n(z;q)\frac{t^{2n}}{[2n+1]_q!},\]
 and 
  \[t\coth_qt \frac{\sinh_q(zt)}{t}=\sum_{n=0}^{\infty}{p}_n(z;q)\frac{t^{2n}}{[2n+1]_q!},\]
  which coincides with the results of Theorem \ref{Sec2:gf}.
\end{exa}
\vskip1mm
\subsection{Relationship with $q$-Appell polynomial sequences}~\par
\vskip1mm

\noindent Recall that the Appell polynomials,  $(a_n(z))_n$ are defined as the polynomials having the series representation 
\begin{equation*}
a_n(z)= \sum_{k=0}^{n} \binom{n}{k}\, a_{n-k}(0)z^k.
\end{equation*}
They are equivalently defined  by the generating function 
\begin{equation*}
g(t)e^{zt}= \sum_{n=0}^{\infty} a_{n}(z)\frac{t^n}{n!},
\end{equation*}
where  $g(t):=\sum_{n=0}^\infty a_n(0) \frac{t^n}{n!}$, $|t|\leq M$, for some $M>0$.   
\vskip3mm

A basic $q$-analog of Appell sequences was first introduced by Sharma and Chak  in~\cite{Sharma-Chak}.  They called them $q$-harmonic. Later, Al-Salam \cite{Al-Salam} studied these sets and referred to it as $q$-Appell polynomials. They defined the $q$-Appell polynomials  $(a_n(z;q))_n$ by
\begin{equation*}
A_q(t)E_q(zt)= \sum_{n=0}^{\infty} a_{n}(z;q)\frac{t^n}{[n]_q!}, \quad n\in \mathbb{N}_0,
\end{equation*}
where  $E_q(z)$ the exponential function defined  in \eqref{E(z)&e(z)},  and $A_q(t)$ is the determining function of the $q$-Appell polynomials given by
\begin{equation*}
A_q(t)=\sum_{n=0}^\infty a_n \frac{t^n}{[n]_q!}.
\end{equation*}
Furthermore, the $q$-Appell polynomials satisfy the following properties:
\begin{enumerate}
  \item $a_0(z;q)\neq 0$;
  \vskip .2 cm
  \item $D_qa_n(z;q)= [n]_q a_{n-1}(z;q)$;
  \vskip .2 cm   
  \item $a_n(0;q)=a_n,\, n\in \mathbb{N}_0$.
\end{enumerate}

\noindent The following result gives characterization of $q$-Appell polynomial sequence (see \cite{Patrick}):
\begin{prop}
A $q$-polynomial sequence $\{a_n(z;q)\}_n$ is a $q$-Appell polynomial sequence if and only if there exists a sequence $(a_k)_{k\geq 0}$; independent of $n$; $a_0\neq 0$, such that
\begin{equation*}
a_n(z;q)= \sum_{k=0}^{n} \qbinom{n}{k}\,q^{\binom{n-k}{2}
} a_{k}z^{n-k}.
\end{equation*}
\end{prop}
Now, we establish a relationship between odd $q$-Lidstone  polynomials of type I and $q$-Appell polynomial sequences.

\begin{thm}\label{$q$-odd Lidstone and q-Appell}
Let $\{a_n(z;q)\}_{n\in \mathbb{N}}$ be a sequence of $q$-Appell polynomials. If $a_{2n+1}(0;q)=0$ $(n\in \mathbb{N})$, then the sequence $(f_n)_n$, where
\begin{equation}\label{relation Appell and Lidstone}
f_n(z;q):= 2^{2n+1} a_{2n+1}(\frac{z}{2};q) \quad (n\in\mathbb{N})
\end{equation}
is an  odd $q$-Lidstone polynomial sequence of type I.
\end{thm}

\begin{proof}
Since $a_{2n+1}(0;q)=0$ and for $q$-Appell polynomial $a_n(z;q)$ the $q$-difference equation
\begin{equation*}
D_qa_n(z;q)= [n]_q a_{n-1}(z;q)
\end{equation*}
 holds,  then  simple calculations yield  that the function $f_n(z;q)$ satisfies \eqref{Def. of qOLS}. This completes the proof of the theorem. 
\end{proof}

\vspace{6pt}

\section{\bf{Even $q$-Lidstone Polynomial Sequences of Type I}}
In this section,  we study the set of polynomial sequences $\{\omega_n(z;q)\}_n$ which satisfy
\begin{equation}\label{Def. of qELS}
\left\{
  \begin{array}{ll}
D_{q}^{2}\omega_n(z;q)= [2n]_q\,[2n-1]_q\, \omega_{n-1}(z;q), & \hbox{} \\\\
D_{q}\omega_n(0;q)=0 \; (n\in \mathbb{N}_0), \quad  \omega_0(z;q)=\gamma_0, \;\gamma_0\in \mathbb{R}\backslash \{0\}. & \hbox{}
  \end{array}
\right.
\end{equation}
We will call these polynomial sets  even $q$-Lidstone polynomial sequences of type I ($q$ELS-I).
\begin{rem}\label{degree qELP}
From  \ref{Def. of qELS}, one can verify that  $\omega_n(z;q)$ is a polynomial of degree $2n$ for all $n\in \mathbb{N}_0$.
\end{rem}

\noindent In this section, The proofs are omitted because they follow along similar lines as the corresponding proofs of Section \ref{Sec:odd}. 
\begin{prop}
The polynomial sequence $\{\omega_n(z;q)\}_n$  is an element of $q$ELS-I if and only if there exists a sequence $\{ \gamma_{2k}\}_{k\geq 0}$ of real numbers such that $\gamma_0 \neq 0$ and
\begin{equation}\label{formala of qELS}
\omega_n(z;q)= \sum_{k=0}^n \qbinom{2n}{2k}\, \gamma_{2(n-k)}\, z^{2k}, \quad  n\in \mathbb{N}.
\end{equation}
\end{prop}
\begin{rem}
From \eqref{formala of qELS}, we obtain
\begin{description}
  \item[i] $\omega_n(z;q)$ is an even function for all $n\in \mathbb{N}$;\\
  \item[ii] $\{z^{2n}\}_n\in qELS$ and  $qELS \subset \hat{P}$, where $\hat{P}:=$ span$\{z^{2j}| \, j\in \mathbb{N}\}$;\\
\item[iii] $ \int_0^1 \omega_n(z;q)\, d_qz=  \sum_{k=0}^n \qbinom{2n}{2k}\, \frac{\gamma_{2(n-k)}}{[2k+1]_q}, \quad  n\in\mathbb{N}$.
\end{description}
\end{rem}
\begin{prop}\label{Prop.1-E}
Let $n\in \mathbb{N}$ and $\{\omega_n(z;q)\}_n\in q$ELS-I. Then
\begin{enumerate}
         \item $D_{q}^{2m} \omega_n(z;q)= \frac{[2n]_q!}{[2(n-m)]!}\, \omega_{n-m}(z;q),\quad  m=0,1,\ldots,n; $\\
         \item $D_{q}^{2m+1} \omega_n(z;q)= \frac{[2n]_q!}{[2(n-m)]!}\, D_{q}\omega_{n-m}(z;q),\,  m=0,1,...,n-1;$\\
         \item $D_{q}^{2m+1} \omega_n(0)=0$ and  $D_{q}^{2m} \omega_n(0)= \frac{[2n]_q!}{[2(n-m)]!}\, \gamma_{2(n-m)},\quad  m=1,2,\ldots,n$.
         \end{enumerate}
\end{prop}

\vskip1mm
\subsection{Matrix form}~\par
\vskip1mm

\begin{defn}\label{matrix form-E}
The $q$-analog of even Lidstone-type I matrix is an infinite lower triangular matrix $F_q=[f_{ij}]$ ($i,j\in \mathbb{N}_0$) with
\begin{equation}\label{elements of A-E}
f_{ij}= \qbinom{2i}{2j}\, \gamma_{2(i-j)}, \quad i\geq j,
\end{equation}
where $\{\gamma_{2k}\}_{k\geq 0}$ a sequence of real numbers and $\gamma_0\neq 0$.
\end{defn}
\begin{rem}
The polynomials \eqref{formala of qELS} can be written in  the matrix form 
\begin{equation}\label{matrix of Elp}
\Omega_q= F_q\hat{Z}_q,
\end{equation}
where $\Omega_q$ and $\hat{Z}$ are two vectors defined by
\begin{equation*}
\Omega_q= [\omega_0(z;q), \omega_1(z;q),\ldots, \omega_n(z;q), \ldots]^T, \quad \hat{Z}=[1, z^2, \ldots, z^{2n}, \ldots]^T.
\end{equation*}
Moreover, if we set $F_{q,n}=[f_{ij}]$ such that $j=0,1,...,i$, $i=0,1, ...,n$ for $n\in \mathbb{N}$ and $f_{ij}$ defined in \eqref{elements of A-E}, then we have a sequence $(F_{q,n})_{n}$ of the principle submatrices of order $n$ of $F_q$ which satisfy
\begin{equation}\label{finite matrix of Elp}
\Omega_{q,n}= F_{q,n}\hat{Z}_{n},
\end{equation}
where
\begin{equation}\label{vector spaces P and Z-E}
\Omega_{q,n}= [\omega_0(z;q), \omega_1(z;q),\ldots, \omega_n(z;q)]^T, \quad \hat{Z}_{n}=[1, z^2, \ldots, z^{2n}]^T.
\end{equation}
\end{rem}
\vskip3mm

\noindent  Let  $T_{q\gamma}$ be a $q$-Toeplitz matrix whose $(i,j)$ entry defined by
\begin{equation*}
t_{ij}^{\gamma}=\left\{
                  \begin{array}{ll}
                    \frac{\gamma_{2(i-j)}}{[2(i-j)]!}, & \hbox{$ i\geq j$;} \\
                    0, & \hbox{otherwise,}
                  \end{array}
                \right.
\end{equation*}
and  $\hat{\mathcal{D}}$  is  the diagonal matrix with entries $\hat{d}_{ii}= [2i]!$.
\begin{prop}
An  even $q$Lidstone-type I matrix, $F_q$,  can be factorized as
\begin{equation}\label{factorize F_q-E}
 F_q= \hat{\mathcal{D}}T_{q\gamma} \hat{\mathcal{D}}^{-1}.
\end{equation}

\end{prop}
\vskip .5 cm 

\begin{prop}\label{(F_q)^{-1}-E}
The $q$-analog of odd Lidstone-type matrix $F_q$ is invertible and
\begin{equation}\label{2factorize F_q-E}
 (F_q)^{-1}= \hat{\mathcal{D}}T_{q\xi} \hat{D}^{-1},
\end{equation}
where $(\xi_{2n})_n$ is the numerical sequence  that satisfies 
\begin{equation}\label{Ebeta sequence}
 \sum_{j=0}^n \frac{\gamma_{2j}\xi_{2(n-j)}}{[2j]_q![2(n-j)]!}= \delta_{n0}, \quad n=0,1, \ldots,
\end{equation}
with $\delta_{nj}$ is the Kronecker's delta.
 \end{prop}

\begin{rem}
Equation \eqref{Ebeta sequence} describes  an infinite linear system which determines the numerical sequence $(\xi_{2n})_n$ and
according to Cramer's rule, the first $n + 1$ equations give
\begin{equation}\label{matrix of Ebeta}
\begin{split}
\xi_0=&\frac{1}{\gamma_0},\\
\xi_{2n}=& (-1)^n\frac{[2]_q![4]_q! \ldots [2n]_q!}{\gamma_0^{n+1}}\\
\times& \mbox{det}\left[
          \begin{array}{ccccc}
             \frac{\gamma_2}{[2]_q!} & \frac{\gamma_0}{[2]_q!} & 0 &  \ldots & 0 \\
             \frac{\gamma_4}{[4]_q!} & \frac{\gamma_2}{[2]_q![2]_q!} & \frac{\gamma_0}{[4]_q!} & 0  & 0 \\
            \vdots & \vdots & \vdots & \ddots & \vdots \\
            \vdots & \vdots & \vdots & \ddots  & \vdots \\
             \frac{\gamma_{2(n-1)}}{[2n-2]_q!} & \frac{\gamma_{2(n-2)}}{[2n-4]_q![2]_q!} & \frac{\gamma_{2(n-3)}}{[2n-6]_q![4]_q!} & \ldots & \frac{\gamma_{0}}{[2n-2]_q!} \\
           \frac{\gamma_{2n}}{[2n]_q!}  & \frac{\gamma_{2(n-1)}}{[2n-2]_q![2]_q!} & \frac{\gamma_{2(n-2)}}{[2n-4]_q![4]_q!} & \ldots & \frac{\gamma_{2}}{[2n-2]_q![2]_q!} \\
          \end{array}
        \right].
\end{split}
\end{equation}
\end{rem}
\vskip5mm

As in the odd $q$-Lidstone polynomial sequences of type I, we consider the $q$-polynomials
\begin{equation}\label{formala of qELS of beta}
\hat{\omega}_n(z;q)= \sum_{k=0}^n \qbinom{2n}{2k}\, \xi_{2(n-k)}\, z^{2k},
\end{equation}
where $(\xi_{2n})_n$ is  the numerical sequence defined in \eqref{matrix of Ebeta}. The two sequences
$\{\omega_n(z;q)\}_n$ and $\{\hat{\omega}_n(z;q)\}_n$  are called conjugate even $q$-Lidstone-type I sequences.

\vskip5mm
\noindent We denote $G_q=[g_{ij}]$,  $( i,j\in \mathbb{N}_0)$ the infinite lower triangular matrix with
\begin{equation*}
g_{ij}=  \qbinom{2i}{2j}\, \xi_{2(i-j)}, \quad i\geq j.
\end{equation*}
Set $\widehat{\Omega}_q= [\hat{\omega}_0(z;q), \hat{\omega}_1(z;q),\ldots, \hat{\omega}_n(z;q), \ldots]^T$. Then, 

\begin{equation}\label{matrix of Elp-beta}
\widehat{\Omega}_q= G_q\hat{Z},
\end{equation}
\begin{equation}\label{matrix of Elp-beta-n}
\widehat{\Omega}_{q,n}= G_{q,n}\hat{Z}_{n}\quad  (n\in \mathbb{N}).
\end{equation}

\begin{prop}\label{conjugate $q$-even Lidstone-type sequences}
The  sequences $\{\omega_n(z;q)\}_n$ and $\{\hat{\omega}_n(z;q)\}_n$  are  conjugate even $q$-Lidstone-type I sequences if and only if
\begin{equation*}
\Omega_q= F_q^2\widehat{\Omega}_q, \quad  \widehat{\Omega}_q= G_q^2\Omega_q,
\end{equation*}
and for $n\in \mathbb{N}$,
\begin{equation*}
\Omega_{q,n}= F_{q,n}^2\widehat{\Omega}_{q,n}, \quad  \widehat{\Omega}_{q,n}= G_{q,n}^2\Omega_{q,n}.
\end{equation*}
\end{prop}
\begin{rem}
From Proposition \ref{conjugate $q$-even Lidstone-type sequences}, we can write
\begin{equation*}
\omega_n(z;q)= \sum_{k=0}^n \tilde{f}_{nk}\, \hat{\omega}_n(z;q), \quad
\hat{\omega}_n(z;q)= \sum_{k=0}^n \tilde{g}_{nk}\, \omega_n(z;q) \, (n\in \mathbb{N}_0),
\end{equation*}
where $\tilde{f}_{nk}$ and $\tilde{g}_{nk}$, $k=0,\ldots, n$, are the elements of the matrices $F^2$ and $G^2$, respectively.
\end{rem}

\vskip1mm
\subsection{Recurrence relations and determinant form}~\par
\vskip1mm

\noindent For the conjugate even $q$-Lidstone-type I polynomial sequences, the following  recurrence relations hold.
\begin{thm}\label{First recu. relation Th.-E}
Let $\{\omega_n(z;q)\}_n$ and $\{\hat{\omega}_n(z;q)\}_n$  ($n\in \mathbb{N}_0$) be the conjugate even $q$-Lidstone-type I polynomial sequences. Then
\begin{equation}
\begin{gathered}
\omega_n(z;q)= \frac{1}{\xi_0}\Big[z^{2n}-\sum_{k=0}^{n-1} \qbinom{2n}{2k}\, \xi_{2(n-k)}\, \omega_k(z;q)\Big], \\
\hat{\omega}_n(z;q)= \frac{1}{\gamma_0}\Big[z^{2n}-\sum_{k=0}^{n-1} \qbinom{2n}{2k}\, \gamma_{2(n-k)}\, \hat{\omega}_k(z;q)\Big],
\end{gathered}
\end{equation}
where $(\gamma_{2k})_n$ and  $(\xi_{2n})_n$ are  the numerical sequences satisfying  \eqref{matrix of Ebeta}.
\end{thm}

\begin{thm}
Let $\{\omega_n(z;q)\}_n\in q$ELS-I. Then
\begin{equation}\label{a first determinate form-E}
\begin{split}
\omega_0(z;q)=&\frac{1}{\xi_0}z,\\
\omega_{n}(z;q)=& \frac{(-1)^n}{\xi_0^{n+1}}
\left|
          \begin{array}{cccccc}
             1 & z^2 & z^4 & \ldots & z^{2n-2} & z^{2n} \\
             \xi_0 & \xi_2 & \xi_4 & \ldots & \xi_{2(n-2)} & \xi_{2n} \\
            0 & \xi_0 &  \qbinom{4}{2}\xi_2 & \ldots & \qbinom{2n-2}{2}\xi_{2(n-2)}& \qbinom{2n}{2}\xi_{2(n-1)} \\
            \vdots & \ddots & \ddots & & \vdots & \vdots \\
\vdots & \ddots & \ddots & & \vdots & \vdots \\
 0 & \ldots & \ldots & & \xi_0 & \qbinom{2n}{2(k-1)}\xi_2
          \end{array}
        \right|,
\end{split}
\end{equation}
where $(\xi_{2n})_n$  is the numerical sequence defined as in \eqref{matrix of Ebeta}. Moreover,  $\{\hat{\omega}_n(z;q)\}_n$ can be expressed in a determinate form similar to \eqref{a first determinate form-E} with $\gamma_{2k}$ instead of $\xi_{2k}$, for $k=0, 1, ..., n$,  and $n\in \mathbb{N}_0$.
\end{thm}
\vskip3mm

\begin{lem}
Let $F_q=[f_{ij}] \, (i,j\in \mathbb{N}_0)$ be an even $q$-Lidstone-type I matrix,  the matrix  $G_q=[g_{ij}]$ be the inverse of $F_q$,  and $\bar{\Pi}_G=[\bar{\pi}_{ij}]$ be the production matrix of $G_q$. Then
\begin{equation}\label{elements of prod.-E}
\begin{split}
\bar{\pi}_{ij}=& \sum_{n=0}^i f_{in} g_{(n+1)j}\\
=&\left\{
    \begin{array}{ll}
      0, & \hbox{$j > i+1$;}  \\
      \sum_{n=0}^{i} \qbinom{2i}{2n} \qbinom{2(n+1)}{2j}\, \gamma_{2(i-n)}\xi_{2(n-j+1)}, & \hbox{otherwise},
    \end{array}
  \right.
\end{split}
\end{equation}
where $(\gamma_{2n})_n$ and $(\xi_{2n})_n$ the numerical sequences defined as in \eqref{Ebeta sequence}.
\end{lem}

\begin{thm}\emph{(Second recurrence relation)}
Let $\{\omega_n(z;q)\}_n$ be a sequence in the class ($q$ELS-I). If $F_q$ is the even $q$-Lidstone-type I matrix related to  $\{\omega_n(z;q)\}_n$,  and
  $\bar{\Pi}_G=[\bar{\pi}_{ij}]\, (i,j\in \mathbb{N}_0)$ is the production matrix of $F^{-1}_q=G_q$, then
\begin{equation}\label{Second recu. relation-E}
\begin{split}
\omega_0(z;q)=&\frac{1}{\xi_0}, \\
\omega_{n+1}(z;q)=& \frac{1}{\bar{\pi}_{n(n+1)}}\Big[z^{2}\omega_{n}(z;q)-\sum_{k=0}^{n} \bar{\pi}_{nk} \omega_k(z;q)\Big], \, n\in \mathbb{N}_0.
\end{split}
\end{equation}
\end{thm}

\begin{thm}
Let $\{\omega_n(z;q)\}_n$ be a sequence in the class ($q$ELS-I).  Then
 \begin{equation*}
\begin{split}
\omega_0(z;q)=&\frac{1}{\xi_0},\\
\omega_{n+1}(z;q)=& \frac{(-1)^{n+1}\omega_0(z;q)}{\bar{\pi}_{01}\bar{\pi}_{12} \ldots \bar{\pi}_{n(n+1)}}\\
\times& \mbox{det}\left[
          \begin{array}{cccccc}
           \bar{\pi}_{00}-z^2 & \bar{\pi}_{01} & 0 & \ldots & \ldots & 0 \\
             \bar{\pi}_{10} & \bar{\pi}_{11}-z^2 & \bar{\pi}_{12} & \ldots & \ldots & 0 \\
            \bar{\pi}_{20} & \bar{\pi}_{21} & \bar{\pi}_{22}-z^2 & \ddots & \ddots & 0 \\
            \vdots & \vdots & \vdots &  \ddots & \ddots & \vdots\\
\vdots & \vdots & \vdots & \ddots & \ddots & \bar{\pi}_{(n-1)n} \\
 \bar{\pi}_{n0} & \bar{\pi}_{n1} & \bar{\pi}_{n2} & \ldots & \ldots & \bar{\pi}_{nn}-z^2
          \end{array}
        \right],
\end{split}
\end{equation*}
where $\bar{\pi}_{ij}$ are  defined as in \eqref{elements of prod.-E}.
\end{thm}

\vskip1mm
\subsection{The Generating function and relationship with $q$-Appell polynomial sequences}~\par
\vskip1mm

\noindent Let $\{\omega_n(z;q)\}_n$ be an even $q$-Lidstone polynomial sequence of type I related to the numerical sequence $(\gamma_{2n})_{n\geq 0}$, and consider the following power series
\begin{equation}\label{power series g-E}
h_q(t)=\sum_{n=0}^\infty \gamma_{2n} \frac{t^{2n}}{[2n]_q!}.
\end{equation}

\begin{lem}\label{Lem.power series 1/g-E}
Let $h_q(t)$ be the power series defined in \eqref{power series g-E}. Then $\frac{1}{h_q(t)}$ is a well-defined function and it has the series representation
\begin{equation}\label{power series 1/g-E}
\frac{1}{h_q(t)}=\sum_{n=0}^\infty \xi_{2n} \frac{t^{2n}}{[2n]_q!},
\end{equation}
where $(\xi_{2n})_n$ is  the numerical sequence defined as in \eqref{matrix of Ebeta}.
\end{lem}

\begin{thm}\label{Sec3:gf}
Let $\{\omega_n(z;q)\}_n$ and $\{\hat{\omega}_n(z;q)\}_n$  be the conjugate even $q$-Lidstone-type I sequences. Then
\begin{eqnarray*}
 h_q(t)\cosh_q (zt)= \sum_{n=0}^\infty \omega_n(z;q) \frac{t^{2n}}{[2n]_q!}, \\
\frac{1}{h_q(t)} cosh_q (zt)= \sum_{n=0}^\infty \hat{\omega}_n(z;q) \frac{t^{2n}}{[2n]_q!},
\end{eqnarray*}
where $h_q(t)$ the function defined as in \eqref{power series g-E}.
\end{thm}

The following result gives a relationship between even $q$-Lidstone and $q$-Appell polynomial sequences.
\begin{thm}\label{$q$-even Lidstone and q-Appell}
Let $(a_n(z;q))_{n\in \mathbb{N}}$ be a sequence of $q$-Appell polynomials. If  $a_{2n+1}(0;q)=0$ $(n\in \mathbb{N})$, then the sequence $(\omega_n)_n$, where
\begin{equation}\label{relation Appell and Lidstone-E}
\omega_n(z;q):= 2^{2n} a_{2n}(\frac{z}{2};q)
\end{equation}
is an even $q$-Lidstone polynomial sequence of type I.
\end{thm}

\begin{exa}
 From \eqref{IM},  
  We  can take 
\[\gamma_{2j}=(-1)^j 2^{2j}\beta_{2j}, \quad \xi_{2j}=\frac{T_{2j+1}}{[2j+1]_q},\]
in \eqref{Ebeta sequence}. 
Then the polynomials $(\widetilde{w}_n(z;q))_n$ and $(\hat{w}_n(z;q))_n$ defined by 
\[w_n(z;q):=\sum_{k=0}^{n}\qbinom{2n}{2k} (-1)^k 2^{2k}\beta_{2k} z^{2n-2k},\]
\[w_n(z;q):=\sum_{k=0}^{n}\qbinom{2n}{2k} \frac{T_{2n-2k+1}}{[2n-2k+1]_q} z^{2k},\]
are conjugate even $q$-Lidstone polynomials of type I. 
Since 
\[\text{Tan}_q t=\tan_qt=\sum_{n=0}^{\infty}T_{2n+1}\frac{t^{2n+1}}{[2n+1]_q!} ,\]
 \[t \text{Cot}_q t=t\cot_qt=\sum_{n=0}^{\infty}(-1)^n\beta_{2n}(q)\frac{ t^{2n}}{[2n]_q!},\]
 we can prove that 
 \[\frac{\tan_qt}{t}{\cosh}_q(zt)=\sum_{n=0}^{\infty}\widetilde{w}_n(z;q)\frac{t^{2n}}{[2n]_q!},\]
 and 
  \[t\coth_qt {\cosh}_q(zt)=\sum_{n=0}^{\infty}\hat{w}_n(z;q)\frac{t^{2n}}{[2n]_q!},\]
  which coincides with the results of Theorem \ref{Sec3:gf}.
\end{exa}
\vskip .5 cm 
\section{\bf{Odd and even $q$-Lidstone Polynomial Sequences of Type II}}

In this section, we  consider two general classes of $q$-polynomial sequences: $\left(\tilde{p}_n(z;q)\right)_n$ and  $\{\tilde{\omega}_n(z;q)\}_n$. We will call these respectively odd and even $q$-Lidstone polynomial sequences of type II. In this type, a sequence of polynomials $\{L_n(z;q)\}_n$ satisfies the $q$-difference equation
\begin{equation*}
D^2_{q^{-1}} L_n(z;q)= a_n\, L_{n-1}(z;q), \quad a_n\in\mathbb{R} \, (n\in \mathbb{N}).
\end{equation*}

In analogy with type I, we give only the statement of some theorems and properties.
\vskip .5 cm 

\subsection{Odd $q$-Lidstone polynomial sequences of type II}~\par
\vskip1mm
\begin{defn}
 The odd $q$-Lidstone sequences of type II ($q$OLS-II) are  the set of polynomial sequences that  satisfy
\begin{equation}
\left\{
  \begin{array}{ll}
D_{q^{-1}}^{2}\tilde{p}_n(z;q)= [2n]_{q^{-1}}\,[2n+1]_{q^{-1}}\, \tilde{p}_{n-1}(z;q), & \hbox{} \\\\
 \tilde{p}_n(0;q)=0 \, (n\in \mathbb{N}_0), \quad  \tilde{p}_0(z;q)=\tilde{\alpha}_0 z, \, \tilde{\alpha}_0\in \mathbb{R}\backslash \{0\}. & \hbox{}
  \end{array}
\right.
\end{equation}
\end{defn}
Notice that   we may assume that the sequence $(\tilde{p}_n)_n$ satisfies the $q$-difference  equation
\[D_{q^{-1}}^{2}\tilde{p}_n(z;q)= [2n]_q\,[2n+1]_q\, \tilde{p}_{n-1}(z;q),\] without  any  loss of generality.

\noindent Similar to Equation \eqref{formala of qOLS}, the following result gives a characterization  of the set $q$OLS-II.
\begin{prop}
The $q$-analog of Lidstone-type sequence $\left(\tilde{p}_n(z;q)\right)_n$ is an element of $q$OLS-II if and only if there exists a sequence $\{ \tilde{\alpha}_{2k}\}_{k\geq 0}$ of real numbers such that $\tilde{\alpha}_0 \neq 0$,  and
\begin{equation}\label{formala of qOLS-II}
\begin{split}
\tilde{p}_n(z;q)=& \sum_{k=0}^n \qbinom{2n+1}{2k+1}\, \frac{q^{2(n-k)^2+(n-k)}\,\tilde{\alpha}_{2k}}{[2(n-k)+1]_q}\,  z^{2(n-k)+1}\\
=& \sum_{k=0}^n \qbinom{2n+1}{2k+1}\, \frac{ q^{k(2k+1)}\, \tilde{\alpha}_{2(n-k)}}{[2(n-k)+1]_q}\, z^{2k+1}.
\end{split}
\end{equation}
\end{prop}
\vskip3mm
 The matrix form  of the odd $q$-Lidstone sequences of type II  is an infinite lower triangular matrix $\tilde{A}_q=[\tilde{a}_{ij}]$ ($i,j\in \mathbb{N}_0$) with
\begin{equation*}
\tilde{a}_{ij}= \qbinom{2i+1}{2j+1}\, \frac{q^{j(2j+1)} \, \tilde{\alpha}_{2(i-j)}}{[2(i-j)+1]_q}, \, i\geq j,
\end{equation*}
where $\{\tilde{\alpha}_{2k}\}_{k\geq 0}$ is  a sequence of real numbers,  and $\tilde{\alpha}_0\neq 0$. Notice,  Formula \eqref{formala of qOLS-II} can be written in the  matrix form
\begin{equation*}
\tilde{P}_q= \tilde{A}_q Z,
\end{equation*}
where $\tilde{P}_q$ and $Z$ are two vectors defined by
\begin{equation*}
\tilde{P}_q= [\tilde{p}_0(z;q), \tilde{p}_1(z;q),\ldots, \tilde{p}_n(z;q), \ldots]^T, \quad Z=[z, z^3, \ldots, z^{2n+1}, \ldots]^T.
\end{equation*}

\begin{prop}
The odd $q$-Lidstone-type II matrix $\tilde{A}_q$ can be factorized as
\begin{equation*}
 \tilde{A}_q= \mathcal{D}T_{q\tilde{\alpha}} \mathcal{D}^{-1},
\end{equation*}
where $T_{q\tilde{\alpha}}=[t_{ij}^{\tilde{\alpha}}]$ with
\begin{equation*}
t_{ij}^{\tilde{\alpha}}= \left\{
                           \begin{array}{ll}
                             \frac{q^{j(2j+1)} \,\tilde{\alpha}_{2(i-j)}}{[2(i-j)+1]_q!}, & \hbox{$i\geq j$;} \\
                             0, & \hbox{ otherwise},
                           \end{array}
                         \right.
\end{equation*}
and $\mathcal{D}$  is the diagonal matrix with entries $d_{ii}= [2i+1]_q!$. Moreover, the matrix $\tilde{A}_q$ is invertible and
\begin{equation*}
 (\tilde{A}_q)^{-1}= \mathcal{D}T_{q\tilde{\beta}} D^{-1},
\end{equation*}
where $(\tilde{\beta}_{2n})_n$ is the numerical sequence satisfying 
\begin{equation}\label{beta sequence II}
 \sum_{j=0}^n \frac{\tilde{\beta}_{2j}\tilde{\alpha}_{2(n-j)}}{[2j+1]_q![2(n-j)+1]_q!}= \delta_{n0}, \quad n\in \mathbb{N}_0.
\end{equation}
\end{prop}

\begin{defn}(Conjugate odd $q$-Lidstone polynomials)
Let $(\tilde{\alpha}_{2n})_n$ and $(\tilde{\beta}_{2n})_n$ be two numerical sequences satisfying  Equation \eqref{beta sequence II}, and $\left(\hat{p}_n(z;q)\right)_n$ be the polynomials
\begin{equation*}
\hat{p}_n(z;q)= \sum_{k=0}^n \qbinom{2n+1}{2k+1}\, \frac{q^{k(2k+1)} \,\tilde{\beta}_{2(n-k)}}{[2(n-k)+1]_q}\, z^{2k+1}, \quad n\in \mathbb{N}_0.
\end{equation*}
Then, the two sequences $\left(\tilde{p}_n(z;q)\right)_n$ and $\left(\hat{p}_n(z;q)\right)_n$ are called conjugate odd $q$-Lidstone sequences of type II.
\end{defn}
\begin{rem}
If $\tilde{B}_q=[\tilde{b}_{ij}]$, $( i,j\in \mathbb{N}_0)$  is the infinite lower triangular matrix with entities 
\begin{equation*}
\tilde{b}_{ij}=  \qbinom{2i+1}{2j+1}\, \frac{q^{j(2j+1)} \,\tilde{\beta}_{2(i-j)}}{[2(i-j)+1]_q}, \quad i\geq j,
\end{equation*}
and $\tilde{P}_q^*= [\hat{p}_0(z;q), \hat{p}_1(z;q),\ldots, \hat{p}_n(z;q), \ldots]^T$, then we have the matrix form
\begin{equation*}
\tilde{P}_q^*= \tilde{B}_qZ.
\end{equation*}
\end{rem}
\begin{prop}
The  sequences $\left(\tilde{p}_n(z;q)\right)_n$ and $\left(\hat{p}_n(z;q)\right)_n$  are  conjugate odd $q$-Lidstone sequences of type II if and only if
\begin{equation*}
\tilde{P}_q= \tilde{A}_q^2\tilde{P}_q^*, \;\mbox{and}\; \tilde{P}_q^*= \tilde{B}_q^2\tilde{P}_q.
\end{equation*}
\end{prop}
\vskip3mm

In the following theorem,  we determine  the recurrence relations for the sequences $\left(\tilde{p}_n(z;q)\right)_n$ and  
$\left(hat{p}_n(z;q)\right)_n$.
\begin{thm}
Let $\left(\tilde{p}_n(z;q)\right)_n$ and $\left(\hat{p}_n(z;q)\right)_n$ be conjugate odd $q$-Lidstone sequences of type II. Then,
\begin{equation*}
\begin{split}
\tilde{p}_n(z;q)=& \frac{q^{-n(2n+1)}}{\tilde{\beta}_0}\Big[z^{2n+1}-\sum_{k=0}^{n-1} \qbinom{2n+1}{2k+1}\,q^{k(2k+1)} \frac{\tilde{\beta}_{2(n-k)}}{[2(n-k)+1]_q}\, \tilde{p}_k(z;q)\Big];\\
\hat{p}_n(z;q)=& \frac{q^{-n(2n+1)}}{\tilde{\alpha}_0}\Big[z^{2n+1}-\sum_{k=0}^{n-1} \qbinom{2n+1}{2k+1}\, q^{k(2k+1)}\frac{\tilde{\alpha}_{2(n-k)}}{[2(n-k)+1]_q}\, \hat{p}_k(z;q)\Big],
\end{split}
\end{equation*}
where $(\tilde{\alpha}_{2n})_n$ and $(\tilde{\beta}_{2n})_n$ are the numerical sequences satisfying  Equation \eqref{beta sequence II}.
\end{thm}
\begin{cor}
The  polynomials $\left(\tilde{p}_n(z;q)\right)_n$ and $\left(\hat{p}_n(z;q)\right)_n$ satisfy the $q$-difference equations
\begin{equation*}
\begin{gathered}
\sum_{k=0}^n \frac{\tilde{\beta}_{2k}}{[2k+1]_q!}\, q^{(n-k+1)(2n-2k+3)} D_{q^{-1}}^{2k}u(z)-  z^{2n+1}=0,\\
\sum_{k=0}^n \frac{\tilde{\alpha}_{2k}}{[2k+1]_q!}\, q^{(n-k+1)(2n-2k+3)} D_{q^{-1}}^{2k}u(z)-  z^{2n+1}=0.
\end{gathered}
\end{equation*}
\end{cor}
\begin{lem}
Let $\tilde{A}_q=[\tilde{a}_{ij}]$ be an odd $q$-Lidstone-type II matrix, $\tilde{B}_q=[\tilde{b}_{ij}]$ be the inverse of $\tilde{A}_q$,  and  $\tilde{\Pi}_B=[\tilde{\pi}_{ij}]$ be the production matrix of $\tilde{B}_q$. Then
\begin{equation}\label{PI-II}
\tilde{\pi}_{ij}= \left\{
    \begin{array}{ll}
      \tilde{\alpha}_0\tilde{\beta}_2, & \hbox{$i=j=0$}, \\
      0, & \hbox{$j > i+1$},  \\
     \sum_{n=0}^{i-j+1} \qbinom{2i+1}{2(n+j)-1}\, q^{(n+j-1)(2(n+j)-1}\frac{\tilde{\beta}_{2n}\tilde{\alpha}_{2(i-j-n)+2}[2(n+j)+1]_q!}{(2(i-j-n)+3)[2j+1]_q![2n+1]_q!}, & \hbox{otherwise},
    \end{array}
  \right.
\end{equation}
where $(\tilde{\alpha}_{2n})_n$ and $(\tilde{\beta}_{2n})_n$  are the numerical sequences defined as in \eqref{beta sequence II}.
\end{lem}

\begin{thm}\emph{(Second recurrence relation)}
Let $\left(\tilde{p}_n(z;q)\right)_n$ be in the class $q$OLS-II.  Assume that $\tilde{A}_q$ is the odd $q$-Lidstone-type II matrix related to  $\left(\tilde{p}_n(z;q)\right)_n$,  and
  $\tilde{\Pi}_q=[\tilde{\pi}_{ij}]$ is the production matrix of $\tilde{A}^{-1}_q$ defined as in \eqref{PI-II}.  Then
\begin{equation}\label{Second recu. relation II}
\begin{split}
\tilde{p}_0(z;q)=&\frac{1}{\tilde{\beta}_0}z, \\
\tilde{p}_{n+1}(z;q)=& \frac{1}{\tilde{\pi}_{n(n+1)}}\Big[z^{2}\tilde{p}_{n}(z;q)-\sum_{k=0}^{n} \tilde{\pi}_{nk}\, \tilde{p}_k(z;q)\Big],\; \left(n\in \mathbb{N}\right).
\end{split}
\end{equation}
Moreover, the conjugate sequence $\left(\hat{p}_n(z;q)\right)_n$ has a relation similar to \eqref{Second recu. relation II} with $\tilde{\alpha}_0$ instead $\tilde{\beta}_0$, and $\tilde{\Pi}_q$ is the production matrix of $\tilde{A}_q$ instead of $\tilde{A}^{-1}_q$.
\end{thm}

\begin{thm}\emph{(Generating functions)}
Let $\left(\tilde{p}_n(z;q)\right)_n$ and $\left(\hat{p}_n(z;q)\right)_n$  be the conjugate odd $q$-Lidstone-type II sequences. Then
\begin{eqnarray}
 \tilde{g}_q(t)\frac{Sinh_q (zt)}{t}= \sum_{n=0}^\infty \tilde{p}_n(z;q) \frac{t^{2n}}{[2n+1]_q!}, \\
\frac{1}{t\tilde{g}_q(t)}\frac{Sinh_q (zt)}{t}= \sum_{n=0}^\infty \hat{p}_n(z;q) \frac{t^{2n}}{[2n+1]_q!},
\end{eqnarray}
where  $\tilde{g}_q(t)$ the power series defined by
\begin{equation*}
\tilde{g}_q(t)=\sum_{n=0}^\infty \tilde{\alpha}_{2n} \frac{t^{2n}}{[2n+1]_q!}.
\end{equation*}
\end{thm}

Now, consider the $q$-Appell polynomials that  satisfy  the $q$-difference equation:
\begin{equation}\label{App-II}
D_{q^{-1}}\tilde{a}_n(z;q)= [n]_q \tilde{a}_{n-1}(z;q).
\end{equation}

A relationship between $q$-Appell polynomial sequences  and odd $q$-Lidstone polynomial sequences of type II is in the following theorem.

\begin{thm}\label{2$q$-odd Lidstone and q-Appell}
Let $\{\tilde{a}_n(z;q)\}_{n\in \mathbb{N}}$ be a sequence of $q$-Appell polynomials of degree $n$. If $\tilde{a}_{2n+1}(0;q)=0$ $(n\in \mathbb{N})$, then the function
\begin{equation}\label{relation Appell and Lidstone II}
f_n(z;q):= 2^{2n+1} \tilde{a}_{2n+1}(\frac{z}{2};q)
\end{equation}
is in the class of odd $q$-Lidstone polynomial sequences of type II.
\end{thm}

\begin{exa}
Ismail and Mansour~\cite[Eq.(3.37)]{Ismail and Mansour} introduced the identity 
\begin{equation}\label{IM2}
\sum_{k=0}^{n}(-1)^k 2^{2k}\frac{\beta_{2k}(q)}{[2k]_q!} \frac{T_{2n-2k+1}(q)}{[2n-2k+1]_q!}=\delta_{n,0},
\end{equation}
where $\delta_{n,0}$ is the Kroncker's delta.  It is worth noting that there was a small typo in \cite[Eq.(3.37)]{Ismail and Mansour},  which we have corrected in \eqref{IM2}. 
  We  can take 
\[\widetilde{\beta}_{2j}=(-1)^j 2^{2j}[2j+1]_q\beta_{2j}, \quad\widetilde{\alpha}_{2j}=T_{2j+1}.\]
Then the polynomials $({P}_n(z;q))_n$ and $(\hat{P}_n(z;q))_n$ defined by 
\[\hat{P}_n(z;q):=\sum_{k=0}^{n}\qbinom{2n+1}{2k+1} q^{k(2k+1)}(-1)^k 2^{2k}\beta_{2k} z^{2k+1},\]
\[P
_n(z;q):=\sum_{k=0}^{n}\qbinom{2n+1}{2k+1} q^{k(2k+1)}\frac{T_{2n-2k+1}}{[2n-2k+1]_q} z^{2k+1},\]
are conjugate odd $q$-Lidstone polynomials of type II. 
Since 
\[\text{Tan}_q t=\tan_qt=\sum_{n=0}^{\infty}T_{2n+1}\frac{t^{2n+1}}{[2n+1]_q!} ,\]
 \[t \text{Cot}_q t=t\cot_qt=\sum_{n=0}^{\infty}(-1)^n\beta_{2n}(q)\frac{ t^{2n}}{[2n]_q!},\]
 we can prove that 
 \[ t\coth_qt\frac{\text{Sinh}_q(zt)}{t}=\sum_{n=0}^{\infty}{P}_n(z;q)\frac{t^{2n}}{[2n]_q!},\]
 and 
  \[\frac{\tan_qt}{t}\frac{\text{Sinh}_q(zt)}{t}=\sum_{n=0}^{\infty}\hat{P}_n(z;q)\frac{t^{2n}}{[2n+1]_q!},\]
  which coincides with the results of Theorem \ref{Sec4:gf}.
\end{exa}
\vskip1mm
\subsection{Even $q$-Lidstone polynomial sequences of type II}~\par
\vskip1mm

\begin{defn}
 An  even $q$-Lidstone sequences of type II ($q$ELS-II) is a  set of polynomial sequences  satisfying
\begin{equation}
\left\{
  \begin{array}{ll}
D_{q^{-1}}^{2}\tilde{\omega}_n(z;q)= [2n]_q\,[2n-1]_q\, \tilde{\omega}_{n-1}(z;q), & \hbox{} \\\\
 \tilde{\omega}_n(0;q)=0 \, (n\in \mathbb{N}_0), \quad  \tilde{\omega}_0(z;q)=\tilde{\gamma}_0 z, \, \tilde{\gamma}_0\in \mathbb{R}\backslash \{0\}. & \hbox{}
  \end{array}
\right.
\end{equation}
\end{defn}

\begin{prop}
A  sequence of polynomials, $\{\tilde{\omega}_n(z;q)\}_n$,  is in the class   $q$ELS-II if and only if there exists a sequence $\{ \tilde{\gamma}_{2k}\}_{k\geq 0}$ of real numbers such that $\tilde{\gamma}_0 \neq 0$, and
\begin{equation}\label{formala of qELS-II-}
\tilde{\omega}_n(z;q)= \sum_{k=0}^n \qbinom{2n}{2k}\, q^{k(2k-1)} \tilde{\gamma}_{2(n-k)}\, z^{2k}.
\end{equation}
\end{prop}
\vskip3mm
 
 Note that  the identity in  \eqref{formala of qELS-II-} can be written in the  matrix form
\begin{equation*}
\tilde{\Omega}_q= \tilde{F}_q \hat{Z},
\end{equation*}
where
\begin{equation*}
\tilde{\Omega}_q= [\tilde{\omega}_0(z;q), \tilde{\omega}_1(z;q),\ldots, \tilde{\omega}_n(z;q), \ldots]^T, \quad \hat{Z}=[1, z^2, \ldots, z^{2n}, \ldots]^T,
\end{equation*}
and 
$\tilde{F}_q=[\tilde{f}_{ij}]$ ($i,j\in \mathbb{N}_0$) with
\begin{equation*}
\tilde{f}_{ij}= \qbinom{2i}{2j}\, q^{j(2j-1)} \frac{ \tilde{\gamma}_{2(i-j)}}{[2(i-j)+1]_q}, \, i\geq j.
\end{equation*}

\begin{prop}
The even $q$-Lidstone-type II matrix $\tilde{F}_q$ can be factorized as
\begin{equation*}
 \tilde{F}_q= \hat{\mathcal{D}}T_{q\tilde{\gamma}} \hat{\mathcal{D}}^{-1},
\end{equation*}
where $T_{q\tilde{\gamma}}=[t_{ij}^{\tilde{\gamma}}]$ with
\begin{equation*}
t_{ij}^{\tilde{\gamma}}= \left\{
                           \begin{array}{ll}
                             \frac{q^{j(2j-1)} \,\tilde{\gamma}_{2(i-j)}}{[2(i-j)]!}, & \hbox{$i\geq j$,} \\
                             0, & \hbox{ otherwise},
                           \end{array}
                         \right.
\end{equation*}
and $\hat{\mathcal{D}}$  is the diagonal matrix with entries $\hat{d}_{ii}= [2i]!$. Moreover, the matrix $\tilde{F}_q$ is invertible and
\begin{equation*}
 (\tilde{F}_q)^{-1}= \hat{\mathcal{D}}T_{q\tilde{\xi}} \hat{D}^{-1},
\end{equation*}
where $(\tilde{\xi}_{2n})_n$  is  the numerical sequence satisfying
\begin{equation}\label{beta sequence II-}
 \sum_{j=0}^n \frac{\tilde{\xi}_{2j}\tilde{\gamma}_{2(n-j)}}{[2j]_q![2(n-j)]!}= \delta_{n0} \quad (n\in \mathbb{N}_0),
\end{equation}
and $\delta_{n0}$ is the Kronecker's delta.
\end{prop}

\begin{defn}(Conjugate even $q$-Lidstone polynomials)
Let $(\tilde{\gamma}_{2n})_n$ and $(\tilde{\xi}_{2n})_n$ be two numerical sequences satisfying  Equation \eqref{beta sequence II-}, and $\{\hat{\omega}_n(z;q)\}_n$ be the polynomials
\begin{equation*}
\hat{\omega}_n(z;q)= \sum_{k=0}^n \qbinom{2n}{2k}\, q^{k(2k-1)} \,\tilde{\xi}_{2(n-k)}\, z^{2k} \quad( n\in \mathbb{N}_0).
\end{equation*}
Then, the two sequences $\{\tilde{\omega}_n(z;q)\}_n$ and $\{\hat{\omega}_n(z;q)\}_n$ are called conjugate even $q$-Lidstone sequences of type II.
\end{defn}
\begin{rem}
If $\tilde{G}_q=[\tilde{g}_{ij}]$,  $( i,j\in \mathbb{N}_0)$ is  the infinite lower triangular matrix with
\begin{equation*}
\tilde{g}_{ij}=  \qbinom{2i}{2j}\, q^{j(2j-1)} \,\tilde{\xi}_{2(i-j)}, \quad i\geq j,
\end{equation*}
and $\tilde{\Omega}^*_q= [\hat{\omega}_0(z;q), \hat{\omega}_1(z;q),\ldots, \hat{\omega}_n(z;q), \ldots]^T$, then we have the  matrix form representation
\begin{equation*}
\tilde{\Omega}^*_q= \tilde{G}_q\hat{Z}.
\end{equation*}
\end{rem}
\begin{prop}
The  sequences $\left(\tilde{\omega}_n(z;q)\right)_n$ and $\left(\hat{\omega}_n(z;q)\right)_n$  are  conjugate even $q$-Lidstone sequences of type II if and only if
\begin{equation*}
\tilde{\Omega}_q= \tilde{F}_q^2\tilde{\Omega}^*_q, \quad  \tilde{\Omega}_q^*= G_q^2\tilde{\Omega}_q.
\end{equation*}
\end{prop}
\vskip3mm

In the following theorem,  we determine  recurrence relations for the sequences $\{\tilde{\omega}_n(z;q)\}_n$,  and  $\{\hat{\omega}_n(z;q)\}_n$.
\begin{thm}
Let $\{\tilde{\omega}_n(z;q)\}_n$ and $\{\hat{\omega}_n(z;q)\}_n$ be conjugate even $q$-Lidstone sequences of type II. Then
\begin{equation*}
\begin{split}
\tilde{\omega}_n(z;q)=& \frac{q^{n(1-2n)}}{\tilde{\xi}_0}\Big[z^{2n}-\sum_{k=0}^{n-1} \qbinom{2n}{2k}\,q^{k(2k-1)} \tilde{\xi}_{2(n-k)}\, \tilde{\omega}_k(z;q)\Big];\\
\hat{\omega}_n(z;q)=& \frac{q^{n(1-2n)}}{\tilde{\gamma}_0}\Big[z^{2n}-\sum_{k=0}^{n-1} \qbinom{2n}{2k}\, q^{k(2k-1)}\tilde{\gamma}_{2(n-k)}\, \hat{\omega}_k(z;q)\Big],
\end{split}
\end{equation*}
where $(\tilde{\gamma}_{2n})_n$ and $(\tilde{\xi}_{2n})_n$ are the numerical sequences defined in Equation \eqref{beta sequence II-}.
\end{thm}

The following theorem  gives the generating functions of the even $q$-Lidstone polynomial sequence of type II.

\begin{thm}\label{Sec4:gf}
Let $\left(\tilde{\omega}_n(z;q)\right)_n$,  and $\left(\hat{\omega}_n(z;q)\right)_n$  be the conjugate even $q$-Lidstone-type II sequences. Then
\begin{eqnarray}
 \tilde{h}_q(t) Cosh_q (zt)= \sum_{n=0}^\infty \tilde{\omega}_n(z;q) \frac{t^{2n}}{[2n]_q!}, \\
\frac{1}{\tilde{h}_q(t)} Cosh_q (zt)= \sum_{n=0}^\infty \hat{\omega}_n(z;q) \frac{t^{2n}}{[2n]_q!},
\end{eqnarray}
where  $\tilde{h}_q(t)$ is  the power series defined by
\begin{equation*}
\tilde{h}_q(t)=\sum_{n=0}^\infty \tilde{\gamma}_{2n} \frac{t^{2n}}{[2n]_q!}.
\end{equation*}
\end{thm}
\begin{exa}
 From \eqref{IM},  
  We  can take 
\[\widetilde{\gamma}_{2j}=(-1)^j 2^{2j}\beta_{2j}, \quad \widetilde{\xi}_{2j}=\frac{T_{2j+1}}{[2j+1]_q},\]
in \eqref{beta sequence II-}. 
Then the polynomials $(\widetilde{w}_n(z;q))_n$ and $(\hat{w}_n(z;q))_n$ defined by 
\[w_n(z;q):=\sum_{k=0}^{n}\qbinom{2n}{2k}q^{k(2k-1)} (-1)^{n-k} 2^{2n-2k}\beta_{2n-2k} z^{2k},\]
\[\hat{w}_n(z;q):=\sum_{k=0}^{n}\qbinom{2n}{2k}q^{k(2k-1)} \frac{T_{2n-2k+1}}{[2n-2k+1]_q} z^{2k},\]
are conjugate even $q$-Lidstone polynomials of type II. 
Since 
\[\begin{gathered}\sum_{n=0}^{\infty}w_n(z;q)\frac{t^{2n}}{[2n]_q!}=\sum_{n=0}^{\infty}\frac{t^{2n}}{[2n]_q!}\sum_{k=0}^{n}\qbinom{2n}{2k}q^{k(2k-1)} (-1)^{n-k} 2^{2n-2k}\beta_{2n-2k} z^{2k}\\
=\left(\sum_{k=0}^{\infty}\frac{T_{2n+1}}{[2n+1]_q!}t^{2n}\right)\left(\sum_{n=0}^{\infty}\frac{q^{n(2n-1)}}{[2n]_q!} (zt)^{2k}\right)\\
=\frac{\tan_q t }{t}\text{Cosh}_q t .\end{gathered}\]
Similarly, we can prove that 
\[
\sum_{n=0}^{\infty}\hat{w}_n(z;q)\frac{t^{2n}}{[2n]_q!}= (t\cot_q t) \text{Cosh}_q t . 
\]

  which coincides with the results of Theorem \ref{Sec4:gf}.
\end{exa}

The even $q$-Lidstone polynomial sequences of type II have a relationship with $q$-Appell polynomials similar to
Theorem \ref{$q$-odd Lidstone and q-Appell}.

\begin{thm}
Let $\left(\tilde{a}_n(z;q)\right)_{n\in \mathbb{N}}$ be a sequence of $q$-Appell polynomials that  satisfy \eqref{App-II}. If  $\tilde{a}_{2n+1}(0;q)=0$ $(n\in \mathbb{N})$, then the function
\begin{equation*}
\tilde{\omega}_n(z;q)= 2^{2n} \tilde{a}_{2n}(\frac{z}{2};q)
\end{equation*}
is in the class  of  even $q$-Lidstone polynomial sequences of type II.
\end{thm}

\section{\bf{Examples}}

Here, we consider some illustrative examples of the odd and even $q$-Lidstone polynomial sequences of type I and type II. These sequences are associated  with $q$-Bernoulli and $q$-Euler's polynomials generated by the first and second Jackson $q$-Bessel functions (see \cite{Ismail and Mansour}).

\begin{exa}\label{Ex.1}
Let $0<q<1$ and  $b_{n}(z;q)$ be the $q$-Bernoulli polynomials defined by the generating function
\begin{equation*}
\dfrac{t\,e_q(zt)}{e_q(t/2)E_q(t/2)-1}=\sum_{n=0}^{\infty}b_n(z;q)
\frac{t^n}{[n]_q!},
\end{equation*}
where $E_q(z)$ and $e_q(z)$ are the $q$-exponential functions defined as in \eqref{E(z)&e(z)}.
Since $D_qb_n(z;q)=[n]_qb_{n-1}(z;q)$, see~\cite{Ismail and Mansour},  then 
 from Theorem \ref{$q$-odd Lidstone and q-Appell},  the sequence $\{p_n\}_n$   defined as
\begin{equation}
p_n(z;q):= 2^{2n+1} b_{2n+1}(\frac{z}{2};q), \quad n\in \mathbb{N}
\end{equation}
is an odd $q$-Lidstone polynomials sequence of type I. Therefore,  it satisfies \eqref{Def. of qOLS}. By Equation \eqref{D_qp_n}, we have $D_qp_n(0;q)=\alpha_{2n}$. This implies
\begin{equation}\label{alpha for B_n}
\alpha_{2n}=  2^{2n} [2n+1]_q \beta_{2n}(q),\quad n\in \mathbb{N}_0,
\end{equation}
where $\beta_{2n}(q)$ denotes the  $q$-Bernoulli numbers. In other words,    $\beta_n(q):=b_n(0;q)$.  Ismail and Mansour in~\cite{Ismail and Mansour} introduced the expansion 
\[t\coth_q t=\sum_{n=0}^{\infty}\beta_{2n}\frac{(2t)^n}{[2n]_q!}.\]
 Consequently, from \eqref{power series g},  
\begin{equation*}
\frac{g_q(t)}{t}= \frac{1}{t}\sum_{n=0}^\infty \alpha_{2n} \frac{t^{2n}}{[2n+1]_q!}= \frac{1}{t}\sum_{n=0}^\infty 2^{2n}\frac{\beta_{2n}(q)}{ [2n]_q!} \,t^{2n}=coth_q(t).
\end{equation*}
Consequently,  from \eqref{g1}, the generating function of the sequence  $\{p_n\}_n$ is $\frac{g_q(t)}{t}\sinh_q zt$, i.e., 
\begin{equation*}
\coth_q(t) \sinh_q(tz)= \sum_{n=0}^\infty p_{n}(z;q) \frac{t^{2n}}{[2n+1]_q!}.
\end{equation*}
Moreover, from \eqref{g2},  the generating function for the conjugate sequence  $\{\hat{p}_n\}_n$ is
\begin{equation*}
\tanh_q(t) \sinh_q(tz)= \sum_{n=0}^\infty \hat{p}_{n}(z;q) \frac{t^{2n+2}}{[2n+1]_q!}.
\end{equation*}
Since 
\[\tanh_q t=\text{Tanh}_qt=-\sum_{n=0}^{\infty}\frac{\widetilde{E}_{2n+1}(q)}{[2n+1]_q!} 2^{2n+1} t^{2n+1},\]
see~\cite[Eq. (3.36) ]{Ismail and Mansour}
Therefore,
\[ \hat{p}_{n}(z;q)=-\sum_{k=0}^{n}\qbinom{2n+1}{2k+1}\frac{\widetilde{E}_{2k+1}(q)}{[2n-2k+1]_q}2^{2k+1} z^{2n-2k+1} \quad (n\in\mathbb{N}_0).\]
\end{exa}

\begin{exa}\label{Ex.2}
Let  $B_{n}(z;q)$ be the $q$-Bernoulli polynomials generated by the second Jackson $q$-Bessel functions which defined in \eqref{B(z;q)}. Consider the sequence $\left(\tilde{p}_n(z;q)\right)_n$ defined by
\begin{equation*}
\tilde{p}_n(z;q)= 2^{2n+1} B_{2n+1}(\frac{z}{2};q).
\end{equation*}
By the same argument as in Example \ref{Ex.1}, the set $\{p_n\}_n$ is an odd $q$-Lidstone polynomials sequence of type II with $
\tilde{\alpha}_{2n}=  2^{2n} [2n+1]_q \beta_{2n}(q)$ for every $n\in \mathbb{N}_0$.
According to Equation \eqref{formala of qOLS-II}, we get
\begin{equation*}
\tilde{p}_n(z;q)= \sum_{k=0}^n \qbinom{2n+1}{2k+1}\, q^{k(2k+1)} 2^{2(n-k)} \beta_{2(n-k)}(q)\, z^{2k+1}.
\end{equation*}
Also,  the conjugate sequences $\left(\tilde{p}_n(z;q)\right)_n$ and $\left(\hat{p}_n(z;q)\right)_n$ have the following generating functions:
\begin{equation*}
\begin{split}
Coth_q(t) Sinh_q(tz)=& \sum_{n=0}^\infty \tilde{p}_n(z;q)\frac{t^{2n}}{[2n+1]_q!}, \\
Tanh_q(t) Sinh_q(tz)=& \sum_{n=0}^\infty \hat{p}_{n}(z;q) \frac{t^{2n}}{[2n+1]_q!}.
\end{split}
\end{equation*}
Thus,  
\[ \hat{p}_{n}(z;q)=-\sum_{k=0}^{n}\qbinom{2n+1}{2k+1}\frac{\widetilde{E}_{2k+1}(q)}{[2n-2k+1]_q}2^{2n-2k+1} q^{k(2k+1)}z^{2k+1} \quad (n\in\mathbb{N}_0).\]
\end{exa}

\begin{rem}
In Example \ref{Ex.2}, the odd $q$-Lidstone polynomial sequence $\{p_n\}_n$ up to a constant $[2n+1]_q!$ coincides with the $q$-Lidstone polynomials $A_n(z)$ defined in \eqref{A_n+B_n}.
\end{rem}

\begin{exa}\label{Ex.Euler}
Let  $E_{n}(z;q)$ be the $q$-Euler polynomials generated by the second Jackson $q$-Bessel functions which defined in \eqref{Euler Second kind}.  Consider the sequence $\{f_n(z;q)\}_n$ defined by
\begin{equation*}
f_n(z;q)= 2^{2n+1} E_{2n+1}(\frac{z}{2};q).
\end{equation*}
Then, the set $\{f_n\}_n$ is an odd $q$-Lidstone polynomials sequence of type II with
\begin{equation*}
\tilde{\alpha}_{2n}=  2^{2n} [2n+1]_q\tilde{E}_{2n}(q),
\end{equation*}
where $\tilde{E}_{2n}(q)= E_{2n}(0;q)$.  Taking into account that $\tilde{E}_{2n}(q)= \delta_{n,0}$ where $\delta_{n,0}$ is the Kronecker's delta, we get
$\tilde{\alpha}_{0}=1$ and $\tilde{\alpha}_{2n}=0$ for every $n\in \mathbb{N}$. This implies $h_q(t)=1$ and then
 the conjugate sequences $\{f_n(z;q)\}_n$ and $\{\hat{f}_n(z;q)\}_n$ have the following generating functions:
\begin{equation*}
\begin{split}
\frac{1}{t} Sinh_q(tz)=& \sum_{n=0}^\infty f_{n}(z;q) \frac{t^{2n}}{[2n+1]_q!}, \\
t \, Sinh_q(tz)=& \sum_{n=0}^\infty \hat{f}_{n}(z;q) \frac{t^{2n}}{[2n+1]_q!}.
\end{split}
\end{equation*}
\end{exa}

\begin{rem}
In Example \ref{Ex.Euler}, the odd $q$-Lidstone polynomials sequence $\{f_n\}_n$ multiplied by  $\frac{2}{[2n+1]_q!}$ coincides with the $q$-Lidstone polynomials $N_n(z)$ defined in \eqref{M_n+N_n}.
\end{rem}

\begin{exa}\label{Ex.3}
 Let  $\{\omega_n(z;q)\}_n$ be the sequence of polynomials defined by
\begin{equation}\label{w_n}
\omega_n(z;q)= 2^{2n} e_{2n}(\frac{z}{2};q), \quad n\in \mathbb{N}_0,
\end{equation}
where  $e_{n}(z;q)$ is the $q$-analog of Euler polynomial defined by the generating function
\begin{equation*}
\dfrac{2\,e_q(zt)}{e_q(t/2)E_q(t/2)+1}=\sum_{n=0}^{\infty}e_n(z;q)
\frac{t^n}{[n]_q!}.
\end{equation*}
By Theorem \ref{$q$-even Lidstone and q-Appell}, the sequence  $\{\omega_n\}_n$ defined in \eqref{w_n} is an even  $q$-Lidstone polynomials sequence of type I, and then it satisfies Equation \eqref{formala of qELS}. According to Proposition \ref{Prop.1-E}, we get
\begin{equation}
\gamma_{2n}= \omega_n(0)= 2^{2n} \tilde{E}_{2n}(q)= 2^{2n}\delta_{n,0}, \quad n\in \mathbb{N}_0.
\end{equation}
 This implies
\begin{equation*}
h_q(t)= \sum_{n=0}^\infty \gamma_{2n} \frac{t^{2n}}{[2n]_q!}= 1.
\end{equation*}
So, the generating function of the sequence  $\{\omega_n\}_n$ is
\begin{equation}\label{gener. of Euler}
cosh_q(tz)= \sum_{n=0}^\infty \omega_{n}(z;q) \frac{t^{2n}}{[2n]_q!}.
\end{equation}
\end{exa}

\begin{exa}\label{Ex.4}
Let  $E_{n}(z;q)$ be the $q$-Euler polynomials generated by the second Jackson $q$-Bessel functions defined by \eqref{Euler Second kind}.
Consider the sequence $\{\tilde{\omega}_n(z;q)\}_n$:
\begin{equation*}
\tilde{\omega}_n(z;q)=  2^{2n} E_{2n}(\frac{z}{2};q), \quad n\in \mathbb{N}_0.
\end{equation*}
By the same argument as in Example \ref{Ex.3}, the set $\{\tilde{\omega}_n(z;q)\}_n$ is an even $q$-Lidstone polynomials sequence of type II with $
\tilde{\gamma}_{2n}= \tilde{E}_{2n}(q)$ for every $n\in \mathbb{N}_0$.
Here, the conjugate sequences $\{\tilde{\omega}_n(z;q)\}_n$ and $\{\hat{\omega}_n(z;q)\}_n$ have the same generating functions:
\begin{equation*}
Cosh_q(tz)= \sum_{n=0}^\infty \tilde{\omega}_n(z;q)\frac{t^{2n}}{[2n]_q!}=\sum_{n=0}^\infty \hat{\omega}_{n}(z;q) \frac{t^{2n}}{[2n]_q!}.
\end{equation*}
\end{exa}

\bibliographystyle{plain}

\bibliographystyle{}

\end{document}